\documentclass[12pt,a4paper]{article}

\usepackage{lineno,hyperref}
\modulolinenumbers[5]
\usepackage{amssymb,amsfonts,amsmath}
\usepackage[all,arc]{xy}
\usepackage{fullpage}
\usepackage{xspace}
\usepackage{amssymb, mathrsfs}
\usepackage{amsmath, graphicx, amsthm}
\usepackage{times}
\usepackage{pgf}
\usepackage{tikz}
\usepackage{tikz-cd}
\usetikzlibrary{matrix,arrows,decorations.pathmorphing}
\usepackage{enumerate}
\usepackage{hyperref}
\hypersetup{
    colorlinks=true,
    linkcolor=blue,
    filecolor=magenta,      
    urlcolor=cyan,
    citecolor=red,
    bookmarks=true}
\usepackage{cleveref}
\usepackage{mathtools}
\usepackage{latexsym}
\usepackage{textcomp}
\usepackage{amscd}
\usepackage{amsxtra}
\usepackage{stmaryrd}
\usepackage{mathrsfs}
\usepackage{todonotes}
\usepackage{ulem}

\newtheorem{thm}{Theorem}[section]
\newtheorem{cor}[thm]{Corollary}
\newtheorem{pr}[thm]{Proposition}

\newtheorem*{Thm}{Theorem}

\newtheorem{sbprop}[subsubsection]{Proposition}
\newtheorem{sbcor}[subsubsection]{Corollary}
\newtheorem{sblem}[subsubsection]{Lemma}

\newenvironment{pf}{\proof[\proofname]}{\endproof}

\theoremstyle{definition}

\newtheorem{para}[subsubsection]{}
\newtheorem{sbdefn}[subsubsection]{Definition}

\newtheorem{sbrem}[subsubsection]{Remark}

\theoremstyle{definition}

\theoremstyle{remark}

\newcommand{\tit}{\textit}

\newcommand{\bt}[1][normal]{\begin{tikzcd}[ampersand replacement = \&, column sep=#1,row sep=#1]}
\newcommand{\etk}{\end{tikzcd}}
\newcommand{\te}[1]{\mathrel{ \text{#1} }}

\newcommand{\cH}{{\mathcal{H}}}

\newcommand{\cO}{{\mathcal{O}}}

\newcommand{\cR}{{\mathcal{R}}}

\newcommand{\sH}{{\mathscr{H}}}

\newcommand{\sM}{{\mathscr{M}}}

\newcommand{\sS}{{\mathscr{S}}}
\newcommand{\sT}{{\mathscr{T}}}

\newcommand{\sV}{{\mathscr{V}}}

\newcommand{\sX}{{\mathscr{X}}}
\newcommand{\sY}{{\mathscr{Y}}}
\newcommand{\sZ}{{\mathscr{Z}}}

\newcommand{\al}{\alpha}

\renewcommand{\th}{\theta}

\newcommand{\la}{\lambda}
\newcommand{\La}{\Lambda}

\newcommand{\si}{\sigma}

\newcommand{\Om}{\Omega}

\newcommand{\R}{\mathbb R}

\newcommand{\Z}{\mathbb Z}

\newcommand{\Q}{\mathbb Q}

\newcommand{\ten}{\otimes}              

\newcommand{\F}{\mathbb{F}}

\newcommand{\js} {\mathcal{J}_{\si}}
\newcommand{\is} {\mathcal{I}_{\si}}

\newcommand{\m}{\mathfrak{m}}

\DeclareMathOperator{\ch}{char}
\DeclareMathOperator{\dl}{dlog} 
\DeclareMathOperator{\rsw}{rsw} 
\DeclareMathOperator{\Gal}{Gal}

\DeclareMathOperator{\Hom}{Hom}

\newcommand{\Bl}{{{\rm Bl}}}

\newcommand{\nl}{{\rm nlog}}

\newcommand{\gp}{{\rm gp}}
\newcommand{\Spec}{{\rm Spec}}

\newcommand{\Spa}{{\rm Spa}}

\makeatletter
\let\c@equation\c@thm
\makeatother
\numberwithin{equation}{section}

\begin{document}

\title{Upper Ramification Groups for Arbitrary Valuation Rings}

\author{Kazuya Kato, Vaidehee Thatte}

\date{}
\maketitle

\begin{abstract} 

T. Saito established a ramification theory for ring extensions locally of complete intersection. We show that for a Henselian valuation ring $A$ with field of fractions $K$ and for a finite Galois extension $L$ of $K$, the integral closure $B$ of $A$ in $L$ is a filtered union of subrings of $B$ which are of complete intersection over $A$. By this, we can obtain a ramification theory of  Henselian valuation rings as the limit of the ramification theory of Saito. Our theory generalizes the ramification theory of complete discrete valuation rings of Abbes-Saito. We study ``defect extensions'' which are not treated in these previous works.

\end{abstract}

\tableofcontents

\section{Introduction}\label{s:in}

\subsection{A brief summary.}\label{Intro1}
Let $A$ be a Henselian valuation ring with field of fractions $K$, let $\bar K$ be a separable closure of $K$, let $G=\Gal(\bar K/K)$, and let $\bar A$ be the integral closure of $A$ in $\bar K$. Since $A$ is Henselian, $\bar A$ is a valuation ring. 

In this paper, for nonzero proper ideals $I$ of $\bar A$, we define closed normal subgroups $G^I_{\log}$ and $G^I_{\nl}$ of $G$  (``$\nl$'' means non-logarithmic or ``non-log'' for short) which we call the upper ramification groups. 
\\\\
We have $G_{\log}^I \subset G_{\nl}^I$, and $G_{\log}^I\supset G_{\log}^J$, and $G_{\nl}^I\supset G^J_{\nl}$ if $I\supset J$ (see \Cref{easy}). 
\\\\
This is a generalization of the work of A. Abbes  and T. Saito (\cite{AS1}) on discrete valuation rings (see \ref{relAS}). 
\\\\
Our work is closely related to the work of T. Saito  \cite{Sa}  (see \ref{withSa}).
\\\\
 A remarkable aspect of this paper which does not appear in the works \cite{AS1} and \cite{Sa} is that the defect can be non-trivial. That is, if the residue field of $A$ is of characteristic $p>0$ and $L$ is a cyclic extension of $K$ of degree $p$, it is possible that the extensions of the residue field  and the value group are both trivial. Our Theorem \ref{thm1} below shows that our upper ramification groups can catch  the defect for such $L/K$.

\subsection{Compatibility for DVRs.}\label{relAS}
For  discrete valuation rings $A$, Abbes-Saito \cite{AS1} defined the logarithmic upper ramification groups $G^r_{\log}$ and their non-log version $G^r$ for $r\in \Q_{>0}$. In the case where the residue field is perfect, if we denote by $G^r_{\text{cl}}$ the classical upper ramification groups, we have $G_{\text{cl}}^r= G^r_{\log}= G^{r+1}$. (For $0 < r \leq 1$,  $G^r$ coincides with the inertia subgroup of $G$.)  \\ 

Their $G^r_{\log}$ (resp. $G^r$) coincides with our $G^{I(r)}_{\log}$ (resp. $G^{I(r)}_{\nl}$) for the ideal $I(r)$ of $\bar A$ defined by $I(r):=\{x\in \bar A\;|\; \text{ord}_{\bar A}(x)\geq r\}$, and our $G_{\log}^I$ (resp. $G^I_{\nl}$) coincides with the closure in $G$ of the union of their $G^r_{\log}$ (resp. $G^r$)  where $r$ ranges over all elements of $\Q_{>0}$ such that $I(r)\subset  I$.

\subsection{Relation with the work \cite{Sa}.}\label{withSa} 
T. Saito  developed the ramification theory  of finite flat rings $B'$ over $A$ which are of complete intersection over $A$ (in \cite{Sa} Section 3.2). 

As we will see in Section \ref{s:re} (Theorem \ref{thm2}), for a finite Galois extension $L$ of $K$ and for the integral closure $B$ of $A$, $B$ is a filtered union of subrings $B'$ of $B$ over $A$ which are  finite flat over $A$ and of complete intersection over $A$. This Theorem \ref{thm2} is deduced from results in \cite{V1,V2,V3}. In this paper, we obtain important results on the upper ramification groups  as the ``limit'' of Saito's ramification theories for $B'/A$.

\subsection{Relation to the ramification theory in \cite{V1,V2,V3}.}\label{withV}
 Assume that the residue field of $A$ is of positive characteristic $p$.
Our upper ramification groups match well with the ramification theories  \cite{V1,V2,V3} of cyclic extensions $L$ of $K$ of degree $p$. In those papers, we considered 
an ideal $\sH$ of $A$, which is a generalization of the classical Swan conductor and plays an important role in the ramification theory of $L/K$.  Some of its crucial properties are :
\begin{itemize} 
\item $\sH=A$ if and only if $L/K$ is unramified. 
\item $\sH$ is not a principal ideal if and only if $L/K$ is a defect extension. 
\item When $K$ is of characteristic $p$, $\sH$ is the ideal generated by all nonzero elements $h$ of $A$ such that $L$ is generated by the solution $\alpha$ of the Artin-Schreier equation  $\alpha^p-\alpha=1/h$. 
\end{itemize}

\noindent (We note that the ideal $\sH$ here was denoted by $H$ and $\cH$ in \cite{V1} and \cite{V2}, respectively.)
\\\\
 We will prove the following result:

\begin{thm} \label{thm1} Assume that the residue field of $A$ is of characteristic $p>0$. Let $L$ be a cyclic extension of $K$ of degree $p$ and let $\sH\subset A$ be the associated ideal. Then for a nonzero proper ideal $I$ of $\bar A$, the image of $G_{\log}^I$ in $\Gal(L/K)$ is $\Gal(L/K)$ if and only if $I\cap A\supset \sH$ and is $\{1\}$ if and only if $I\cap A \subsetneq \sH$.

\end{thm}

\noindent Thus, our upper ramification groups can catch the important ideal $\sH$.

\subsection{Indexing by ideals of $\bar A$.}\label{why} As the index set of the upper ramification filtration, we use the set of all nonzero proper ideals of $\bar A$,  not the positive part of (the value group of $A$) $\otimes \Q$. (The latter is identified with the set of all principal nonzero proper ideals of $\bar A$.) Here we explain the reason.
\\\\
Consider a Henselian valuation ring $A$  whose residue field is of characteristic $p>0$,  cyclic extensions $L_1$ and $L_2$ of degree $p$ of $K$, a nonzero element $a$ of ${\mathfrak m}_{\bar A}$, such that the ideal $\sH$ (\ref{withV}) of $A$ associated to $L_1/K$ generates the ideal  $J_1=a \bar A$ of $\bar A$ and the ideal $\sH$ of  $A$ associated to $L_2/K$ generates the ideal $J_2=a{\mathfrak m}_{\bar A}$ of $\bar A$. (Such $A$, $L_1/K$, $L_2/K$ exist. See \ref{br10}.) By Theorem \ref{thm1}, for a nonzero proper ideal $I$ of $\bar A$,  $G^I_{\log}\to \Gal(L_i/K)$ is surjective if and only if $I\supset J_i$. But if $I$ is  principal, for both $i=1,2$,  $G^I_{\log}\to \Gal(L_i/K)$ is surjective  if and only if $I\supset a{\bar A}$ and thus principal ideals $I$  cannot catch the difference of the important ideals $\sH$.

\subsection{Methodology.}\label{method}

 \subsubsection{Concerning $G^I_{\nl}$.} Our definition of $G^I_{\nl}$  follows the methods of Abbes-Saito in \cite{AS1} and Saito in \cite{Sa} except for the following point. In \cite{AS1}, Abbes and Saito used rigid analytic spaces for their definitions of $G^r_{\log}$ and $G^r$.  In \cite{Sa}, Saito used a scheme theoretic algebraic method. 
We use  adic spaces, which are generalizations of rigid analytic spaces. (Actually, we only use the  ``algebraic part'' of the theory of adic spaces, not the analytic part, as is explained in \ref{adican}. We never use the completed coordinate rings of adic spaces. It might be better to say that we use Zariski-Riemann spaces, rather than adic spaces.) 
 
 We could use the scheme theoretic method in \cite{Sa}. But we prefer our method of using adic spaces because an adic space is a space of valuation rings, and we think that  it is natural to use a space of valuation rings to understand the ramification theory of valuation rings. Also, the open covering \ref{opencov} which appears in our method  connects  principal ideals of $\bar A$ and non-principal ideals of $\bar A$ in a nice way.

\subsubsection{Concerning $G^I_{\log}$.} Our method to define   $G^I_{\log}$ is to  modify the definition of $G^I_{\nl}$ by going to log smooth extensions $K'$ of $K$ and replacing a finite Galois extension $L/K$ by the extensions $LK'/K'$. 
The usefulness of log smooth extensions in the case $A$ is a discrete valuation ring is indicated in  \cite[1.2.1 and 1.2.6]{Sa0}. 
 See Section \ref{s:lo} for the definition of log smooth extensions. Tame finite extensions are regarded as log \'etale extensions and log smooth extensions are more general.

\subsection{Outline.}\label{outline}

We define our  upper ramification groups in Section \ref{s:up}. 
Sections  \ref{s:lo} and \ref{s:ad} are preparations for it. In Section \ref{s:lo}, we consider log smooth extensions of Henselian valuation rings.\\
In Section \ref{s:ad}, we consider adic spaces (Zariski-Riemann spaces). 
\\
 In Section \ref{s:Sa},  we review results of Saito in   \cite{Sa}   which we use in this paper.
\\ 
 In Section \ref{s:re}, we deduce the aforementioned \Cref{thm2} and other general results on extensions of valuation rings from the works \cite{V1, V2,V3}.
 \\
In Section \ref{s:pr}, we 
 prove  properties of our upper ramification groups by using Sections \ref{s:Sa} and \ref{s:re}.
 \\ 
In Section \ref{s:AS}, we consider the relation with  Abbes-Saito theory \cite{AS1} described in  \ref{relAS}.
\\ 
 In Section \ref{s:th1}, we prove Theorem \ref{thm1}. 
\\
In Section \ref{s:br}, we give results on breaks of the logarithmic upper ramification filtration.

\subsection{Notation.}\label{N1}  Let $A$ be a Henselian valuation ring which is not a field, and let $K$ be the field of fractions of $A$. Let $\Gamma_A$ denote the value group $K^\times/A^\times$ of $A$, ${\mathfrak m}_A$  the maximal ideal of $A$, and $k$  the residue field $A/\m_A$ of $A$. Consider $\bar K$, a separable closure of $K$, and let $\bar A$ denote the integral closure of $A$ in $\bar K$. 

For a finite extension $L$ of $K$, we denote the integral closure of $A$ in $L$ by $B$. Note that by the assumption $A$ is Henselian, $B$ is a Henselian valuation ring. The index $[\Gamma_B:\Gamma_A]$ is called the ramification index and is denoted by $e(L/K)$. The degree of the extension of the residue fields, called the inertia degree, is denoted by $f(L/K)$.

\section{Log smooth extensions of Henselian valuation rings}\label{s:lo}

In this section, we consider the notion of a {\it log smooth extension} of Henselian valuation rings (\ref{defls}). We will use this later to define the logarithmic upper ramification groups.\\

\subsection{Preliminaries on commutative monoids.}

\subsubsection{Valuative monoids.}\label{val1} Let $\La$ be an abelian  group, whose group law is written multiplicatively,   and 
let $\sV$ be a submonoid of $\La$. We say $\sV$ is a {\it valuative monoid for $\La$} if for each $x\in \La$, we have either $x\in \sV$ or $x^{-1}\in \sV$. 
\\\\
 Let $\La$ be an abelian group and let $\sM_1$ and $\sM_2$ be submonoids of $\La$. We say $\sM_2$ {\it dominates} $\sM_1$ if $\sM_1\subset \sM_2$ and $\sM_2^\times \cap \sM_1= \sM_1^\times$. 

\noindent Here $(-)^\times$ denotes the group of invertible elements. 

\begin{sblem}\label{val2}\cite{Og}[I.2.4.1] Let $\La$ be an abelian group and let $\sM$ be a  submonoid of $\La$. Then there is a valuative monoid for $\La$ which dominates $\sM$.
\end{sblem}

\begin{sblem}\label{val3} Let $\La$ be an abelian group, let $\La_0$ be a subgroup of $\La$, let $\sV$ be a valuative monoid for $\La$, let $\sV_0$ be a valuative monoid for $\La_0$, and assume that $\sV$ dominates $\sV_0$. Then $\sV\cap \La_0=\sV_0$.

\end{sblem}

\begin{pf} Suppose that $x\in \sV \cap \La_0$. If $x\notin \sV_0$, $x^{-1}$ belongs to $\sV_0$. Hence, $x^{-1}\in \sV^\times \cap \sV_0=\sV_0^\times$. This contradicts $x\notin \sV_0$.

\end{pf}

\subsubsection{Associated valuation rings.}\label{Rp}

Assume that we are given a pair $(\La, \sV)$, where $\La$ is an abelian group which contains $K^\times$ as a subgroup such that $\La/K^\times$ is a free abelian group of finite rank, and $\sV$ is a valuative monoid  for $\La$ which dominates the valuative monoid $\sV_0:=A\smallsetminus \{0\}$ for $K^\times$. 

Let $R:=A\otimes_{\Z[\sV_0]} \Z[\sV]$, where $\Z[-]$ denotes the semi-group ring over the ring $\Z$ of integers. Let $\mathfrak p$ be the ideal of $R$ generated by the image of $\sV\smallsetminus \sV^\times$ in $R$.

\begin{sbprop}\label{Rp2} Let the notation be as in \ref{Rp}. Then $\mathfrak p$ is a prime ideal of $R$. The local ring $R_{\mathfrak p}$ is a valuation ring whose value group is canonically isomorphic to $\La/\sV^\times$ and whose residue field is isomorphic to a rational function field over $k$ in $n$ variables where $n$ is the rank of the free abelian group $\sV^\times/A^\times$.

\end{sbprop}

\begin{pf}
We have $R/\mathfrak p= k\otimes_{\Z[A^\times]} \Z[\sV^\times]$. Since $
\sV^\times/A^\times \overset{\subset}\to  \La/K^\times$, $\sV^\times/A^\times $ is a free abelian group of finite rank. Let $U_1, \dots, U_n$ be elements of $\sV^\times$ whose images in $\sV^\times/A^\times$ form a basis. Then $R/\mathfrak p$ is isomorphic to the Laurent polynomial ring $k[U_1^{\pm 1}, \dots, U^{\pm}_n]$ in 
$n$ variables $U_1,\dots, U_n$. Hence, $\mathfrak p$ is a prime ideal of $R$.  Moreover, the residue field of $\mathfrak p$ is the rational function field $k(U_1, \dots, U_n)$ over $k$ in $n$ variables.

Take a subgroup $\La_1$ of $\La$ such that $\La$ is the direct product of $K^\times$ and $\La_1$. Then $\La_1$ is a free abelian group of finite rank. The homomorphism $\Z[\sV_0]\to \Z[\sV]$ is flat by  \cite[Proposition 4.1]{K1}, because the homomorphism $\sV_0\to \sV$ satisfies the condition (iv) in  Proposition 4.1. 

By this, the map $R=A \otimes_{\Z[\sV_0]} \Z[\sV] \to K\otimes_{\Z[\sV_0]}  \Z[\sV]= K\otimes_{\Z[\sV_0^{\gp}]} \Z[\sV_0^{\gp}\sV]\subset K[\La_1]$ is injective, and $R$ (resp. $\mathfrak p$) is identified with the subset of the Laurent polynomial ring $K[\La_1]$ consisting of all elements $\sum_{\la\in \La_1} c_{\la} \la$ ($c_{\la}\in K$) such that $c_\la\la \in \sV$ (resp. $c_{\la}\la \in \sV\smallsetminus \sV^\times$) for all $\la\in \La_1$
 such that $c_\la\neq 0$. We have the valuation $K[\La_1] \smallsetminus \{0\}\to \La/\sV^\times$ given by  $\sum_{\la} c_{\la}\la \mapsto \text{class}(c_{\mu} \mu)$, where $\mu$ is an element of $\La_1$ such that $c_{\mu}\neq 0$ and $c_{\la}\la c_{\mu}^{-1}\mu^{-1}\in \sV$ for all $\la\in \La_1$, and the inverse image of $\sV/\sV^\times$ under this valuation is $R\smallsetminus \{0\}$.

 We will now prove that the valuation ring of this valuation on $\text{frac}(R)$ coincides with the local ring $R_{\mathfrak p}$ of $R$ at $\mathfrak p$.

Each non-zero element of $R$ is a product of an element of $R \backslash \mathfrak{p}$ and an element of $\sV$. Hence, each non-zero element of $\text{frac}(R)$ is written as $ab^{-1}\la$ with $a, b \in R \backslash \mathfrak{p}$ and $\la \in \La$. This element belongs to the valuation ring if and only if $\la \in \sV$, and hence, if and only if it belongs to $R_{\mathfrak{p}}$.\end{pf}

\subsection{Log smoothness - definitions.}\label{defls}

Let $A'\supset A$ be a Henselian  valuation ring which  dominates $A$.

\begin{sbdefn} We say $A'$ is a {\it log smooth extension of $A$ of rational type} if $A'$ is isomorphic over $A$ to the Henselization of $R_{\mathfrak p}$ for the $R$ and $\mathfrak p$ associated to some $\La$, $\sV$ as in \ref{Rp}. \\(The phrase ``rational type" comes from the fact that the residue field of $A'$ is a rational function field over $k$ and the field of fractions of the ring $R$ in \ref{Rp} is a rational function field over $K$.)
\end{sbdefn}

\begin{sbdefn}\label{tfe} We say $A'$ is a  {\it tame finite extension of $A$}  if the field of fractions $K'$ of $A'$ is a  tamely ramified finite extension of  $K$.\\ That is, a finite extension $K'/K$ such that $[K':K]=e(K'/K)f(K'/K)$, $e(K'/K)$ is invertible in $A$, and the residue field $k'$ of $A'$ is a separable extension of $k$. \\ (Since $A'$ is not necessarily finitely generated as an $A$-module, saying $A'$ is a tame finite extension of  $A$ is abuse of terminology.)
\end{sbdefn}

\begin{sbdefn} 

We say $A'$ is a {\it log smooth extension of $A$} if there is a sequence of Henselian valuation rings $A=A_0\subset A_1\subset \dots \subset A_n=A'$ where each extension $A_i/A_{i-1}$ ($1\leq i\leq n$) is either 
 a log smooth extension of rational type or a  tame finite extension. 
\end{sbdefn}

\noindent We note that when the extension of fraction fields is finite, log smoothness is equivalent to tameness.

\subsection{Composition of log smooth extensions.}
\begin{sblem}\label{lse composition}

(1) If $A'$ is a log smooth extension of $A$ and $A''$ is a log smooth extension of $A'$, then  $A''$ is a log smooth extension of $A$. 
\\
(2) Let $A'$ be a log smooth extension of $A$ and let $K'$ (resp. $k'$) be the field of fractions (resp. residue field) of $A'$. Let $\Gamma_{A'}$ denote the value group of $A'$. Then $\Gamma_{A'}/\Gamma_A$ is a finitely generated abelian group and $$\text{trdeg}(K'/K)= \text{trdeg}(k'/k)+\text{rank}(\Gamma_{A'}/\Gamma_A).$$ 

\noindent Here $\text{trdeg}$ denotes the transcendence degree.

\end{sblem}
\begin{pf} (1) This can be seen by combining the two sequences  $A=A_0\subset A_1\subset \dots \subset A_n=A'$ and  $A'=A'_0\subset A'_1\subset \dots \subset A'_{n'}=A''$ of Henselian valuation rings as described in the definition above.
\\\\
(2) We may assume that $A'/A$ is either a tame finite extension or a log smooth extension of rational type. In the former case,  $\Gamma_{A'}/\Gamma_A$ is finite and  $\text{trdeg}(K'/K)= \text{trdeg}(k'/k)=0$. In the latter case, with the notation of Proposition 2.6, $\text{trdeg}(K'/K)$ is equal to the rank of $\La/K^\times$, 
$\text{trdeg}(k'/k)$ is equal to the rank of $\sV^\times/A^\times$, and we have an exact sequence
$$0\to \sV^\times/A^\times\to \La/K^\times\to \Gamma_{A'}/\Gamma_A\to 0$$
because  $\Gamma_{A'}/\Gamma_A\cong (\La/\sV^\times)/(K^\times/A^\times)\cong \La/(K^\times \sV^\times)$.
\end{pf}

 \subsection{Basic types.}\label{types}

 We describe some simple types of log smooth extensions below. 
 Later we will see that the study of log smooth extensions essentially boils down to understanding these types. They will serve as building blocks for the general case.\\

\noindent {\bf Type 1.}  Tame finite extensions (\Cref{tfe}) are log smooth. \\\\
Both  type 2 and type 3 (whose precise descriptions are given below) are log smooth extensions of rational type. In type 2, the residue extension is of transcendence degree 1 and the extension of the value group is finite. In type 3, the residue extension is trivial. 
\\\\
{\bf Type 2.} Take an  integer $e\geq 1$ and  $a\in K^\times$. Let $\mathfrak p$ be the ideal of $A[U]$ generated by $m_A$. Then $\mathfrak p$ is a prime ideal. The integral closure of the local ring $A[U]_{\mathfrak p}$ in $K(U)((aU)^{1/e})$  is a valuation ring. Let $A'$ be the Henselization of this valuation ring. This $A'$ is a log smooth extension of $A$. In fact, in \ref{Rp}, let  $\La= K^\times\times \La_1$ where $\La_1$ is a free abelian group of rank $1$ with generator $\th$. Consider the valuative monoid $\sV$ for $\La$ consisting of 
$c \th^i$ where $c\in K^\times$, $i\in \Z$ such that $c^ea^i\in A$. Then by identifying $\th^ea^{-1}$ with $U$, $A'/A$ is identified with 
 the associated log smooth extension of rational type. 
 \\
 We call $A'/A$ the log smooth extension of type 2 associated to $(e,a)$. 
 The quotient group $\Gamma_{A'}/\Gamma_A$ is isomorphic to a quotient of $\Z/e\Z$ and is generated by the class of $\th$. 

There are special cases or sub-types of type 2. The second case is considered only when  the residue field of $A$ is of positive characteristic $p$.

{\bf Type 2.1.} The case $e=1$. In this case, $\Gamma_A=\Gamma_{A'}$ and the residue field of $A'$ is $k(U)$. 

{\bf Type 2.2.} The case where $e>1$ is a power of $p$, and the class of $a$ in $\Gamma_A$ is not a $p$-th power. 

Then the residue field of $A'$ is $k(U)$ and $\Gamma_{A'}/\Gamma_A\cong \Z/e\Z$. 
\\\\
{\bf Type 3.} Let  $\La_1$ be a free abelian group of finite rank and take a valuative monoid $\sV_1$ for the product group $\Gamma':= \Gamma_A \times  \La_1$ which contains $(A\smallsetminus \{0\})/A^\times\subset \Gamma_A\subset \Gamma'$ such that $\sV_1^\times=\{1\}$. In \ref{Rp}, let $\La= K^\times \times \La_1$ and let $\sV$ be the inverse image of $\sV_1$ in $\La$ and let $A'/A$ be the associated log smooth extension of rational type. 
Then the  value group of $A'$ is identified with $\Gamma'$ and the residue field of $A'$ coincides with that of $A$.   

Conversely, if $A'/A$ is a log smooth extension of Henselian valuation rings of rational type such that the quotient $\Gamma_{A'}/\Gamma_A$ is torsion free and such that the residue field of $A'$ is that of $A$, then $A'/A$ is of this type 3.

\begin{sblem}\label{R23}  Let $A'/A$ be a log smooth extension of rational type. Then there are extensions  $A=A_0\subset A_1\subset \dots \subset A_n \subset A_{n+1}=A'$  such that  $A_i/A_{i-1}$   for $1\leq i\leq n$ is a log smooth extension of type 2 in \ref{types} and  $A_{n+1}/A_n$  is a log smooth extension of type 3 in \ref{types}.

\end{sblem}

\begin{pf} 

Let $\tilde \Theta \subset \La$ be the inverse image of the torsion part $\Theta$ of $\Gamma_{A'}/\Gamma_A$. We prove \ref{R23} by induction on the rank $n$ of the finitely generated free abelian group $\tilde \Theta/K^\times$.

 If $n=0$, $A'/A$ is of type 3 (\ref{types}). Assume $n\geq 1$.
 
 Take compatible  isomorphisms $\Theta\cong  \oplus_{i=1}^n \Z/e(i)\Z$ and $\tilde \Theta/K^\times \cong \oplus_{i=1}^n \Z$. Let $\vartheta$ be an element of $\tilde \Theta$ whose image in  $\oplus_{i=1}^n \Z $ is the first basis element. 
Write ${\vartheta}^{e(1)}=aU$ with $a\in K^\times$ and $U\in \sV^\times$. Let $A_1/A$ be the log smooth extension of type 2  associated to $(e(1), a)$ (\ref{types}).\\ 
Then we have a unique local homomorphism  $A_1\to A'$ of local $A$-algebras sending $\theta$ to $\vartheta$. Let $K_1$ be the field of fractions of $A_1$ and let $\La_1$ be the push out of $K_1^\times \leftarrow K^\times {\vartheta}^{\Z} \to \La$ in the category of abelian groups and let $\sV_1$ be the image of $A_1^\times \sV$ in $\La_1$. Then $\sV_1$ is a valuative monoid for $\La_1$ which dominates $A_1\smallsetminus \{0\}$, and $A'/A_1$ is identified with the log smooth extension of rational type associated to $(\La_1, \sV_1)$. \\ We  have exact sequences 
$$0\to \Z \to \La/K^\times\to \La_1/K_1^\times\to 0 \te{and} 0\to \Z\to \Gamma_{A'}/\Gamma_A \to \Gamma_{A'}/\Gamma_{A_1}\to 0$$ where the second arrows of these sequences  send $1\in \Z$ to the classes of $\vartheta$. \\ Hence, for the inverse image $\tilde \Theta_1$ of the torsion part $\Theta_1$ of $\Gamma_{A'}/\Gamma_{A_1}$ in $\La_1$, $\tilde \Theta_1/K_1^\times$ is of rank $n-1$. This proves \ref{R23} by induction on $n$. 
\end{pf}

\begin{sblem}\label{221} Let $A'/A$ be a log smooth extension of type 2 in \ref{types}.\\ Then we have  $A=A_0\subset A_1\subset A_2=A'$   where  $A_1/A_0$ is a log smooth extension of either 
of type 2.1 or 2.2 in \ref{types} and $A_2/A_1$ is a  tame finite extension.

\end{sblem}

\begin{pf} Let $e'$ be the largest divisor of $e$ which is invertible in $A$. Let $A_1/A$ be the log smooth extension of type 2 associated to $(e/e', a)$ (\ref{types}). This $A_1$ has the desired properties. 
\end{pf}

\subsection{Log smooth extensions are defectless.}
We will discuss some preliminaries before proving the main result (\ref{lsdl}) of this subsection.
\subsubsection{Defect}\label{defect}  We recall the notion of defect in valuation theory (Cf. \cite{Ku}). 
\\
Let $A'\supset A$ be  a Henselian valuation ring which dominates $A$.
Let $K'$ be the field of fractions of $A'$ and let $k'$ be the residue field of $A'$. 
\\\\
(0) Assume $K'$  is a finite extension of $K$. Then $[K':K]\geq e(K'/K)f(K'/K)$. We say the extension {\it $A'/A$ (or $K'/K$) has no defect (i.e. it is defectless)} if the equality holds in this inequality. 
In fact, Ostrowski's Lemma states that there exists a positive integer $d(K'/K)$ such that  $[K':K] = d(K'/K)  e(K'/K)f(K'/K)$. 
Furthermore, $d(K'/K)$ is $1$ in residue characteristic $0$ and a non-negative integral power of $p$ in residue characteristic $p>0$.

This $d(K'/K)$ (sometimes also denoted by $d(A'/A)$) is known as {\it the defect} of the extension $K'/K$ (or the extension $A'/A$). An extension with non-trivial defect is called a {\it defect extension}.
\\\\
(1) If, in particular, $A'$ is finitely generated as an $A$-module, then $K'/K$ is defectless. This can be shown as follows. 

By \cite[Theorem 7.10 ]{M},  $A'$ is a free $A$-module of rank $[K':K]$. 
Since $A'$ is a valuation ring and $A'/\m_AA'$ is an Artinian ring, $A'/\m_AA'$ is a truncated discrete valuation ring of length $r=[K':K]/[k':k]$.
For an element $x \in A'$ whose image in $A'/\m_AA'$ is a prime element, its powers $x^i; 0 \leq i \leq r-1$ have different classes in $\Gamma_{A'}/\Gamma_A$.
Hence, $[\Gamma_{A'}:\Gamma_A] \geq r$.
\\\\
(2)  Assume  $K'/K$ is of finite transcendence degree $\text{trdeg}(K'/K)$. Then $\text{trdeg}(k'/k)$ and $ \dim_{\Q}(\Q\otimes \Gamma_{A'}/\Gamma_A)$ are finite and we have
the Abhyankar inequality: $$\text{trdeg}(K'/K)\geq \text{trdeg}(k'/k)+ \dim_{\Q}(\Q\otimes \Gamma_{A'}/\Gamma_A).$$ We say {\it $A'/A$ has no transcendence defect} if the equality holds in this inequality.  
\\\\
(3)  Assume $K'/K$ has finite transcendence degree.  Then there exists a finite subset $\sT$ of $(K')^\times$ called a {\it valuation transcendence basis} such that $\sT= \sT_1\coprod \sT_2$, where the classes of the members of $\sT_1$ form a $\Q$-basis in $\Q\otimes \Gamma_{A'}/\Gamma_A$, $\sT_2\subset A'$, and the residue classes of the members of $\sT_2$ form a transcendence basis of $k'$ over $k$. 

If  $A'/A$ has no  transcendence defect in the sense of (2), then any valuation transcendence basis is a transcendence basis
of $K'/K$.
\\\\
(4)  Assume  $K'/K$ is of finite transcendence degree. We say $A'/A$ has {\it no defect (i.e. it is defectless)} if $A'/A$ has no transcendence defect in the sense of (2) and if for all subfields $K_1$ and $K_2$ of $K'$ such that $K\subset K_1\subset K_2\subset K'$ and such that $K_2/K_1$ is a finite extension, the extension $(K_2 \cap A')^h/(K_1 \cap A')^h$  has no defect in the sense of (0). 
Here the superscript $(\cdot)^h$ denotes  Henselization.
In such a case we also say that $K_2/K_1$ is defectless.
\\\\
The statements (5) and (6) follow from    \cite[ Theorem 5.4]{Ku} (details in \ref{56}).
\\\\
(5) {\it Theorem.} 
  Assume that $K'/K$ has finite transcendence degree and that $A'/A$ has no  transcendence defect in the sense of (2). Let $\sT$ be 
 a valuation transcendence basis  of $K'/K$ and let $A_1$ be the Henselization of $A'\cap K(\sT)$. Then the extension $A'/A$ has no defect in the sense of (4) if and only if the extension $A'/A_1$ has no defect in the sense of (4). 
\\\\
(6) {\it Theorem.} Let $A''\supset A'$ be a Henselian valuation ring which dominates $A'$, let $K''$ be the field of fractions of $A''$, and assume that $K''/K$ is of finite transcendence degree. 
Then $A''/A$ has no defect in the sense of (4) if and only if $A'/A$ and $A''/A'$ have no defect in the sense of (4). 
\subsubsection{Proofs of (5) and (6).}\label{56}
We will use the relation between log smoothness and the valuation transcendence basis below, as well as the following claim, to prove (5) and (6).
\\
\\
 Let the notation be as in (3). Let $A_1 \subset A'$ be as in (5). Then $A_1/A$ is a log smooth extension of rational type with $\Lambda= \Z^{\sT} \times K^{\times} \subset K'^{\times}$ and with $\sV=\Lambda \cap A'$. 
\\
If $A'/A$ is a log smooth extension of rational type with $\Lambda$ and $\sV$, then we can choose the valuation transcendence basis $\sT=\sT_1 \coprod \sT_2$ of $K'/K$, where $\sT_1$ and  $\sT_2$ are subsets of $\Lambda$  forming bases of $\Lambda/\sV^{\times}$ and $\sV^{\times}/K^{\times}$, respectively, and $\La'$ is the Henselization of $\La' \cap K(\sT)$ for this $\sT$.

\noindent {\bf Claim 1}: Let the notation be as in (3) and let $\sS \subset \sT$. Let $K_1, K_2$ be subfields of $K'$ such that $K \subset K_1 \subset K_2 \subset K'$ and such that $K_2/K_1$ is a finite extension. Then $K_2(\sS)/K_1(\sS)$ is defectless if and only if $K_2/K_1$ is defectless.
\begin{proof}[Proof of Claim 1] 
Let $\sS_1=\sS \cap \sT_1$ and $\sS_2=\sS \cap \sT_2$. 
For $i=1,2$, let $K_i^h$ (resp. $K_i^h(\sT)$) $\subset K'$  denote the fraction field of the Henselization of $K_i$ (resp. $K_i(\sT)$).

By the explicit construction of the log smooth extension $K_i(\sS)^h/K_i^h$ in \cref{Rp2} for $i=1,2$, the residue field of $K_i(\sS)^h$ is the rational function field over the residue field of $K_i^h$ with variables $\sS_2$, and the value group of $K_i(\sS)^h$ is the product of the value group of $K_i^h$ and the free abelian group with basis $\sS_1$.
\\
Therefore, 
we have
$$e(K_2(\sS)^h/K_1(\sS)^h)=e(K_2^h/K_1^h) \te{and} f(K_2(\sS)^h/K_1(\sS)^h)=f(K_2^h/K_1^h).$$
\\
We also have 

$$K_2(\sS)^h=K_1(\sS)^h\otimes_{K_1^h} K_2^h \te{and} [K_2(\sS)^h:K_1(\sS)^h] = [K_2^h:K_1^h].$$
\\
This proves the claim. \end{proof}

\noindent We will now discuss the proofs of (5) and (6). \\\\ As a consequence of  \cite[Theorem 5.4]{Ku}, for a valuation transcendence basis $\sT$ of $K'/K$, whether $K'/K(\sT)$ is defectless or not is independent of the choice of $\sT$. This fact is used in the proofs below.

\begin{proof}[Proof of (5):]
Let $\sT$ be a valuation transcendence basis of $K'/K$. Assume that $K'/K(\sT)$ is defectless. Let $K \subset K_1 \subset K_2 \subset K'$ such that $K_2/K_1$ is finite and $K_i \cap A'$ is Henselian for $i=1,2$. It is enough to prove that $K_2/K_1$ is defectless.

Take a valuation transcendence basis $\{t_1, \cdots, t_m\}$ of $K_1/K$ and a valuation transcendence basis $\{t_{m+1}, \cdots, t_n\}$ of $K'/K_2$. Then $\{t_1, \cdots, t_n\}$ is a valuation transcendence basis of $K'/K$. By the above independence, $K'/K(t_1, \cdots, t_n)$ is defectless.
Since 
$$K(t_1, \cdots, t_n)\subset K_1(t_{m+1}, \cdots, t_n) \subset K_2(t_{m+1}, \cdots, t_n) \subset K',$$
$K_2(t_{m+1}, \cdots, t_n)/K_1(t_{m+1}, \cdots, t_n)$ is defectless. Hence, $K_2/K_1$ is defectless.
\end{proof}
\begin{proof}[Proof of (6):]
We first prove the `if' part. Take a valuation transcendence basis $\{t_1, \cdots, t_m\}$ of $K'/K$ and a valuation transcendence basis $\{t_{m+1}, \cdots, t_n\}$ of $K''/K'$. Since $K'/K(t_1, \cdots, t_m)$ has no defect, $K'(t_{m+1}, \cdots, t_n)/K(t_1, \cdots, t_n)$ has no defect.
Since $K''/K'(t_{m+1}, \cdots, t_n)$ also has no defect, $K''/K'(t_1, \cdots, t_n)$ has no defect. The rest follows from (5).

Next, we prove the `only if' part. Let  $\{t_1, \cdots, t_m\}$  be a valuation transcendence basis of $K'/K$ and let $\{t_{m+1}, \cdots, t_n\}$  be a valuation transcendence basis of $K''/K'$. Then $\{t_1, \cdots, t_n\}$ is a valuation transcendence basis of $K''/K$. By the assumption, $K''/K(t_1, \cdots, t_n)$ is defectless. Hence, $K''/K'(t_{m+1}, \cdots, t_n)$ is defectless and $K'(t_{m+1}, \cdots, t_n)/K(t_1, \cdots, t_n)$ is defectless. By the former and by (5), $K''/K'$ is defectless. By the latter, $K'/K(t_1, \cdots, t_m)$ is defectless. By (5), $K'/K$ is defectless.
\end{proof}

\begin{sbprop}

\label{lsdl} Let  $A'/A$  be a log smooth extension of Henselian valuation rings. Then $A'/A$ has no defect  in the sense of \ref{defect} (4). 

\end{sbprop}

\begin{pf} By the theorem in \ref{defect} (6) and by Lemmas \ref{R23} and \ref{221}, it is sufficient to prove this in the case $A'/A$ is either one of the types 1, 2.1, 2.2 and 3 of \ref{types}. In the case 1, set $\sT=\emptyset$. In the cases 2.1 and 2.2, set $\sT=\{U\}$. In the case 3, let $\sT$ be a lifting of a basis of $\La_1$ to $(K')^\times$. Then $A'/A_1$ in (5) of \ref{defect} has no defect. This follows from the observations that $A'=A_1$ for the types 2.1 and 3, while $e_{A'/A_1} =[K': \text{frac}(A_1)]$ for  the types 1 and 2.2. 

Hence, by the theorem in \ref{defect} (5), $A'/A$ has no defect in the sense of \ref{defect} (4). 
\end{pf}

\subsection{Preliminaries on differential and log-differential modules.}
\begin{sbdefn}\label{om} \textbf{Differential $1$-Forms $\Omega_{\bullet}^1$. }

\begin{enumerate}[(i)]
\item Let $R$ be a commutative ring. The $R$-module $\Om_R^1$ of
\tit{differential $1$-forms over $R$} is defined as follows:
$\Om_R^1$ is generated by
\begin{itemize}
\item The set $\{ db \mid b \in R \}$ of generators.
\item The relations are the usual rules of differentiation:
For all $b, c \in R$,
\begin{enumerate}
\item (Additivity) $d(b+c)=db+dc$
\item (Leibniz rule) $d(bc)=cdb+bdc$
\end{enumerate}
\end{itemize}
\item For a commutative ring $A$ and a commutative $A$-algebra
$B$, the $B$-module $\Omega_{B/A}^1$ of \tit{relative differential
$1$-forms over $A$} is defined to be the cokernel of the map $B
\ten_A \Om^1_A \to \Om^1_B.$
\end{enumerate}
\end{sbdefn}

\begin{sbdefn}\label{oml} \textbf{Logarithmic Differential $1$-Forms $\Omega_{\bullet}^1(\log)$.}
\begin{enumerate}[(i)]
\item For a valuation ring $A$ with the field of fractions $K$,
we define the \textit{$A$-module $\Omega_{A}^1(\log)$ of logarithmic
differential $1$-forms} as follows: $\Omega_{A}^1(\log)$ is generated by
\begin{itemize}
\item The set $\{ db \mid b \in A \} \cup \{ \dl x \mid x
\in K^{\times}\}$ of generators.
\item The relations are the usual rules of differentiation and
an additional rule: For all $b,c \in A$ and for all $x,y \in
K^{\times},$
\begin{enumerate}
\item (Additivity) $d(b+c)=db+dc$
\item (Leibniz rule) $d(bc)=cdb+bdc$
\item (Log 1) $\dl(xy)=\dl x+\dl y$
\item (Log 2) $b \dl b = db$ for all $0 \neq b \in A$
\end{enumerate}
\end{itemize}
\item Let $L/K$ be an extension of Henselian valued fields, $B$
the integral closure of $A$ in $L$, and hence, a valuation ring.
We define the \textit{$B$-module $\Omega_{B/A}^1(\log)$ of logarithmic relative
differential $1$-forms over $A$} to be the cokernel of the map
$B \ten_A \Omega_{A}^1(\log) \to \Omega_{B}^1(\log).$

\end{enumerate}
\end{sbdefn}
\subsubsection{Exact sequence.}\label{ses}
For $A,K$ as in \ref{oml}, we have the   exact sequence of $k$-modules:
$$\begin{matrix}  0 & \to &  \Omega^1_k & \to & k \otimes_A \Omega^1_A(\log) & \to & k \otimes_\Z \Gamma_A & \to 0\end{matrix}.$$

The first map $\Omega^1_k  \to  k \otimes_A \Omega^1_A(\log)$ is given by $dx \mapsto d\tilde{x}$,  for a lift $\tilde{x}$ of $x$ to $A$. For any two lifts $\tilde{x_1}$ and $\tilde{x_2}$, their difference lies in the maximal ideal of $A$. To prove that this map is well-defined, it is enough to show that $d\tilde{x}=0$ for any $\tilde{x} \in \m_A$. We may assume $\tilde{x} \neq 0$.
In this case, $d\tilde{x}=\tilde{x} \otimes \dl \tilde{x}$. This is $0$ since  $\tilde{x}=0$ in $k$.

We note that $k \otimes_A \Omega^1_A(\log)$ is generated by $\dl x ; x \in K^{\times}$ as a $k$-module.\\ The second map $k \otimes_A \Omega^1_A(\log) \to  k \otimes_\Z \Gamma_A$ is given by $\dl x \mapsto$
class of $x$ in $\Gamma_A$.
Let $x \in A$. If $x \in \m_A$, $x \dl x = dx = 0$, and thus, maps to $0$. 
If $x \in A^{\times}$, $dx \mapsto$ residue class of $x ~ \otimes$ class of $x$.

\begin{sblem}\label{injlog}
Let the notation be as in \ref{oml}(ii). Assume that $L/K$ is a finite extension. Let $l$ be the residue field of $B$. Then the map $l \otimes_A \Omega^1_A (\log) \to l \otimes_B\Omega^1_B (\log) $ is injective if and only if  $B/A$ is tamely ramified.
\end{sblem}

\begin{pf}

We have the commutative diagram of exact sequences

    \begin{center}
    \begin{tikzcd}[row sep=huge]
0 \arrow[r] 
& l\otimes_k \Omega^1_k \arrow[r] \arrow[d, "\alpha" ] & l \otimes_A \Omega^1_A(\log) \arrow[r] \arrow[d, , "\beta"]
& l \otimes_\Z \Gamma_A \arrow[d,  "\gamma"] \arrow[r] & 0 \\
0 \arrow[r]
&  \Omega^1_{l} \arrow[r] &  l \otimes_B\Omega^1_{B}(\log) \arrow[r] 
& l \otimes_\Z \Gamma_{B} \arrow[r] & 0. 
\end{tikzcd}
\end{center}

\noindent The extension $l/k$ is separable if and only if $\alpha$ is injective. This is equivalent to  $\alpha$ being bijective.
The map $\gamma$ is injective if and only if $[\Gamma_B:\Gamma_A]$ is invertible in $A$. This concludes the proof.
\end{pf}
\noindent
\subsection{Conditions which are equivalent to log smoothness.}

\begin{sbprop}\label{eq} Let $A'$ be a Henselian valuation ring over $A$ which dominates $A$. Then the  following conditions are equivalent.
\begin{enumerate}[(i)]
\item $A'$ is a log smooth extension of $A$.

\item There are extensions $A=A_0\subset A_1\subset A_2=A'$ of Henselian valuation rings such that $A_1/A_0$ is a log smooth extension of rational type and $A_2/A_1$ is a tame  finite extension.

\item  The following four properties are satisfied. 
\begin{enumerate}[(a)]
\item $\Gamma_{A'}/\Gamma_A$ is  a finitely generated abelian group. 
\item The residue field $k'$ of $A'$ is a finitely generated field over $k$. 
 
\item The field of fractions $K'$ of $A'$ has finite transcendental degree over $K$ and the extension $A'/A$ has no defect in the sense of (4) in \ref{defect}.  

\item The map $k'\otimes_A \Omega^1_A(\log) \to k' \otimes_{A'} \Omega^1_{A'}(\log)$ is injective.
\end{enumerate}
\end{enumerate}
\end{sbprop}

\subsection{Proof of \ref{eq}.}\label{pfeq} 

The implication 
(ii) $\Rightarrow$ (i)  is clear. 
 We will now prove (i) $\Rightarrow$  (iii) and (iii) $\Rightarrow$   (ii).

\begin{proof}[{\bf Proof of (i) $\Rightarrow$ (iii).}]
It is enough to prove this in the cases where $A'/A$ is of type 1 (tame finite), type 2.1, type 2.2, or type 3.

Properties $(a)$ - $(c)$ follow from  \cref{lse composition} and \cref{lsdl}. 
It remains to show that the map $k'\otimes_A \Omega^1_A(\log) \to k' \otimes_{A'} \Omega^1_{A'}(\log)$ is injective.
We have  the following commutative diagram of exact sequences:

\begin{center}
    \begin{tikzcd}[row sep=huge]
0 \arrow[r] 
& k'\otimes_k \Omega^1_k \arrow[r] \arrow[d, "\alpha" ] & k' \otimes_A \Omega^1_A(\log) \arrow[r] \arrow[d, , "\beta"]
& k' \otimes_\Z \Gamma_A \arrow[d,  "\gamma"] \arrow[r] & 0 \\
0 \arrow[r]
&  \Omega^1_{k'} \arrow[r] &  k' \otimes_{A'} \Omega^1_{A'}(\log) \arrow[r] 
& k' \otimes_\Z \Gamma_{A'} \arrow[r] & 0 
\end{tikzcd}
\end{center}

\noindent We want to prove that the middle vertical map $\beta$ is injective when $A'/A$ is of type 1 (tame finite), type 2.1, type 2.2, or type 3.
\\

The first map $\alpha$ is injective, since $k'/k$ is separable (by definition for type 1, $k'=k(U)$ for types 2.1 and 2.2, and $k'=k$ for type 3). 

We observe that the third map $\gamma$ is also injective for type 2.1 as $\Gamma_A=\Gamma_{A'}$. It is also true for type 1 ($[\Gamma_{A'}:\Gamma_{A}]$ invertible in $A$) and type 3 ($\Gamma_{A'}/\Gamma_{A}$ is torsion-free). Therefore, injectivity of $\alpha$ and $\gamma$ in these cases forces $\beta$ to be injective.

Now it remains to prove the injectivity $\beta$ for type 2.2. 
By the injectivity of $\alpha$ and by the snake lemma, it is sufficient to prove that the $k'$-homomorphism $\text{ker}(\gamma) \to \text{coker}(\alpha)$ is injective. 
We have $\text{ker}(\gamma) \cong k'$ and it is generated by the class of $1 \otimes a$. Since $aU$ is a $p$-th power in $K'^{\times}$, the image of the last element in $\text{coker}(\alpha)= \Omega_{k'/k}$ is $- \dl U$ which is nonzero. This completes the proof.
\end{proof}
\begin{proof}[{\bf Proof of (iii) $\Rightarrow$ (ii).}]

We again consider the  commutative diagram of exact sequences:

\begin{center}
    \begin{tikzcd}
0 \arrow[r] 
& k'\otimes_k \Omega^1_k \arrow[r] \arrow[d] & k' \otimes_A \Omega^1_A(\log) \arrow[r] \arrow[d, ]
& k' \otimes_\Z \Gamma_A \arrow[d] \arrow[r] & 0 \\
0 \arrow[r]
&  \Omega^1_{k'} \arrow[r] &  k' \otimes_{A'} \Omega^1_{A'}(\log) \arrow[r] 
& k' \otimes_\Z \Gamma_{A'} \arrow[r] & 0 
\end{tikzcd}
\end{center}

\noindent We define a subgroup $\La$ of $(K')^\times$. It is the subgroup generated by the following elements $T_i$ ($1\leq i\leq \ell$), $U_i$ ($\ell+1\leq i\leq m$), $T_i$ ($m+1\leq i\leq n$), where $0\leq \ell \leq m \leq n$.\\ Note: When $\ell=0$, we only use $U_i$ ($1\leq i\leq m$) and  $T_i$ ($m+1\leq i\leq n$).
\\\\
We first define $T_i$ ($1\leq i \leq \ell$) and $U_i$ ($\ell+1\leq i \leq m$) in the two cases $\ch k = p$ and $\ch k = 0$ as follows. 
\\
{\bf Case $\ch k = p$ :} Let $\Theta$ be  the $p$-power torsion part of $\Gamma_{A'}/\Gamma_A$ which is isomorphic to $\oplus_{i=1}^{\ell} \Z/p^{e(i)}\Z$ with $\ell\geq 0$ and $e(i)\geq 1$. Fix such an isomorphism. Let $T_i$ ($1\leq i \leq \ell$) be an element of $(K')^\times$ whose image in $\Gamma_{A'}/\Gamma_A$ is the 
element of $\Theta$ corresponding to the $i$-th basis element of $\oplus_{i=1}^{\ell} \Z/p^{e(i)}\Z$. Write $T_i^{p^{e(i)}}=a_iU_i$ where 
$a_i\in K^\times$ and $U_i\in (A')^\times$. Then the classes of $a_i$ ($1\leq i\leq \ell$) form an $\F_p$-basis of the kernel of $\Gamma_A/\Gamma_A^p\to \Gamma_{A'}/(\Gamma_{A'})^p$.  
By the above commutative diagram of exact sequences, we have that  $d\log(U_i)$ are linearly independent in $\Omega^1_{k'/k}$, and hence, the residue classes of $U_i$ in $k'$ form a part of a $p$-basis of $k'$ over $k$. Let 
 $U_i$ ($\ell+1\leq i \leq m$, $m\geq \ell$) be elements of $(A')^\times$ such that the residue classes of $U_i$ ($1\leq i\leq m$) form a $p$-basis of $k'$ over $k$. 
 \\
 {\bf Case $\ch k = 0$ :} Take $\ell=0$. Let $(U_i)_{1\leq i\leq m}$ be a family of elements of $(A')^\times$ whose residue classes form a 
 transcendence basis of $k'$ over $k$. 
 \\\\
 Next  (both in the case $\ch k=p$ and in the case $\ch k=0$), we define $T_i$ ($m+1\leq i \leq n$). Take elements $T_i$ ($m+1 \leq i \leq n$, $n\geq m$) of $(K')^\times$ whose classes in $(\Gamma_{A'}/\Gamma_A)/(\text{torsion part})$ 
 form a basis of this free abelian group. 
 \\\\
Let $\La$ be the subgroup of $(K')^\times$ generated by $K^\times$, $T_i$ ($1\leq i \leq \ell$) (recall that  $\ell=0$ 
 in the case $\ch k=0$), $U_i$ ($\ell+1\leq i\leq m$), $T_i$ ($\ell+1\leq i\leq n$). Let $\sV= \La \cap A'$. Let $A_1$ be the log smooth extension of $A$ of rational type associated to $(\La, \sV)$. 
  By the assumption $A'/A$ has no transcendental defect, $K'$ is an algebraic extension of the field of fractions $K_1$ of $A_1$. Let $k_1$ be the residue field of $A_1$. Then $k'$ is a separable finite extension of  $k_1$ and $\Gamma_{A'}/\Gamma_{A_1}$ is a finite group whose order is invertible in $A$. By the assumption $A'/A$ has no defect, this shows that $A'/A_1$ is a tame finite extension.

\end{proof}

\subsection{Further results.}

\begin{sblem}\label{A'''} 
 Let $A_1/A$ and $A_2/A$ be log smooth extensions.  
Then there is $A'$ which is a log smooth extension  of both of $A_1$ and $A_2$. 

\end{sblem}

\begin{pf}  By \ref{eq} (i) $\Rightarrow$ (ii), it is sufficient to prove \ref{A'''} in the case $A_i$ ($i=1,2$) are log smooth extensions of rational type. Assume $A_i$ is associated with  $(\La_i, \sV_i)$  ($i=0,1,2$ and $A=A_0$)  as in \ref{Rp}.
Let $\La$ be the push out of $\La_1\supset K^\times \subset \La_2$. Let $\sM:=\sV_1\sV_2\subset \La$. \\

\noindent {\bf Claim 2}. $\sM$ dominates $\sV_i$ for $i=1, 2$.

\noindent {\it Proof of Claim 2:}  We consider the case of $\sV_1$. Let $x\in \sM^\times\cap \sV_1$. Then $xv_1v_2=1$ for some $v_i\in \sV_i$. Going to $\La/\La_1$, we see that the image of $v_2$ in $\La/\La_1$ is $1$. Hence, the image of $v_2$ in $\La_2/\La_0\cong \La/\La_1$ is $1$, and $v_2\in \sV_1\cap \La_0=\sV_0$ by Lemma \ref{val3}. This shows  $v_1v_2\in \sV_1$, and hence, $x\in \sV_1^\times$.\\

\noindent By Lemma \ref{val2}, there is a valuative monoid $\sV$ for $\La$ which dominates $\sM$. By Claim 2, $\sV$ dominates $\sV_1$ and $\sV_2$. Let $A'$ be the log smooth extension of $A$ of rational type associated to $(\La, \sV)$. 
\end{pf}

\begin{sbprop}\label{124}  Let $A'/A$   be a log smooth extension of Henselian valuation rings, let $K'$ be the field of fractions of $A'$, let $L$ be a finite separable extension of $K$,  let $L'=LK'$, let $B$ be the integral closure of $A$ in $L$, and let $B'$ be the integral closure of $A'$ in $L'$. 
Assume that $\Gamma_A=\Gamma_B$. Then $B'$ is generated by $B$ as an $A'$-algebra and $\Gamma_{A'}=\Gamma_{B'}$. Furthermore, if the residue field of $B$ is  purely inseparable over the residue field of $A$, then the canonical map $B\otimes_A A'\to B'$ is an isomorphism.

\end{sbprop}

\begin{pf} Assume first that $A'$ is a log smooth extension of $A$ of rational type associated to $(\La, \sV)$. Let $\La_B$ be the pushout of $\La\leftarrow K^\times \to L^\times$ in the category of abelian groups and let $\sV_B$ be the inverse image of $\sV/A^\times$ under $\La_B \to \La_B/B^\times \cong \La/A^\times$. Then $B'$ is the log smooth extension of $B$ associated to $(\La_B, \sV_B)$ and from this we have $B\otimes_A A'\overset{\cong}\to B'$. 

Assume next that $K'/K$ is a finite tame extension. Replacing $K$ by the maximal unramified extension of $K$ in $K'$, we may assume that $[\Gamma_{A'}:\Gamma_A]=[K':K]=n$. We have a basis $(e_i)_{1\leq i\leq n}$ of the $K$-vector space $K'$ such that the valuations of $e_i$ form a set of representatives of $\Gamma_{A'}/\Gamma_A$ in $\Gamma_{A'}$. Let $x\in B'$. Write $x=\sum_{i=1}^n x_ie_i$ where $x_i\in L$. Then $x_ie_i\in B'$ for all $i$.  Write $x_i=u_iy_i$ where $u_i\in B^\times$ and $y_i\in K$. Then $y_ie_i \in A'$. This shows that $B\otimes _A A' \to B'$ is an isomorphism.   \end{pf}

\begin{sblem}\label{tame'}
 Let $A'/A$   be a log smooth extension of Henselian valuation rings, let $K'$ be the field of fractions of $A'$, let $L$ be a finite extension of $K$. If $LK'/K'$ is tamely ramified then $L/K$ is also tamely ramified.
\end{sblem}

\begin{pf}
    Let $B, B'$ be the integral closures of $A, A'$ in $L$ and  $LK'$, respectively. Let $l'$ be the residue field of $B'$.
We have the following commutative diagram:

    \begin{center}
    \begin{tikzcd}[row sep=huge]
 l' \otimes_A \Omega^1_A(\log) \arrow[r] \arrow[d]
& l' \otimes_B \Omega^1_B(\log) \arrow[d]  \\
l' \otimes_{A'} \Omega^1_{A'}(\log) \arrow[r] 
& l' \otimes_{B'} \Omega^1_{B'}(\log) 
\end{tikzcd}
\end{center}

\noindent The vertical arrows in the diagram are injective by log smoothness. The lower horizontal map is injective by tameness (\cref{injlog}). Hence, the upper horizontal map is injective. Applying \cref{injlog} again, this proves that $L/K$ is tamely ramified.
\end{pf}

\section{Adic spaces (Zariski-Riemann spaces)}\label{s:ad}
 Let $A$ and $K$ be as in \ref{N1}. Let $L$ be a finite Galois extension of $K$ and let $B$ be the integral closure of $A$ in $L$. 

In the case $A$ is a complete discrete valuation ring, to study the ramification in the extension $L/K$, Abbes and Saito used rigid analytic spaces  $X_Z^a$ (\cite{AS1},  3.1). 
For a general  $A$, we use locally ringed spaces $X(S, T, I)$ defined as below, 
which  play the roles of $X_Z^a$. 
Here $(S, T)$ is a pre-presentation of $B/A$ in the sense of \ref{presen} below and  $I$ is  a non-zero proper ideal of $\bar A$.

\subsection{Presentation and pre-presentation.}\label{presen}

Let $S$ be a 
 finite subset of $B$. Let $T$ is a finite subset of the kernel of the homomorphism $A[\{y_s\}_{s\in S}]\to B\;;\; y_s\mapsto s$, where $y_s$ are indeterminates.

\begin{sbdefn}\label{KP}

 Pre-presentation. We call such a  pair $(S, T)$  a  {\it pre-presentation of $B/A$}. We call it also a {\it pre-presentation for $L/K$}.\end{sbdefn}

\begin{sbdefn}\label{pres} Presentation. We call $(S, T)$ a {\it presentation of $B/A$, or a presentation for $L/K$},
if the map $A[\{y_s\}_{s\in S}]/(T) \to B$ is an isomorphism.
\end{sbdefn}

\subsection{The unit disk $D^n$.}

 We explain the definition of $X(S,T,I)$ briefly but not precisely. The precise definition will be given in \ref{XSTI}.\\\\ If $n=\sharp(S)$,  $X(S, T, I)$ is an open subspace of the $n$-dimensional unit disc $D^n$ over $\bar A$ (defined below)
  consisting of all points of $D^n$ at which we have $t\equiv 0\bmod I$ for all $t\in T$. 
 \\\\
In the case $A$ is a complete discrete valuation ring, Abbes and Saito use a presentation $(Z, T)$ of $B/A$, and  their $X_Z^a$ for $a\in \Q_{> 0}$  coincides with our  $X(Z, T, I)$  as a topological space, where $I$ is the ideal of $\bar A$ generated by an $a$-th power of a prime element of $A$.  
\\\\
For a general  $A$, a presentation of $B/A$ need not exist  (the ring $B$ need not be  finitely generated  over $A$), and therefore, we use a pre-presentation. 
\\\\
We describe the $n$ dimensional unit disc $D^n$ over $\bar A$ in \ref{sZ}--\ref{cODn}.

\subsubsection{}\label{sZ}  Let $$R= {\bar K}[T_1, \dots, T_n]\supset R^+= {\bar A}[T_1, \dots, T_n].$$

\noindent  Let $Z$ be the projective limit of the blowing ups $\Bl_J(X)$ of $X=\Spec(R^+)$ along all finitely generated ideals $J$ of $R^+$ such that $RI=R$.\\\\
 Here, the projective system is made by $\Bl_{JJ'}(X) \to \Bl_J(X)$ for such ideals $J$ and $J'$. We regard $Z$ as a locally ringed space as follows. The topology of $Z$ is the projective limit of the Zariski topologies of $\Bl_J(X)$,  and the structure sheaf $\cO_Z$ on $Z$ is the inductive limit of the inverse images of the structure sheaves of $\Bl_J(X)$.
\subsubsection{}
Note that the  pullbacks of the blowing ups $\Bl_J(X) \to X$ in \ref{sZ} to  $\Spec(R)$  
are isomorphisms, and hence, the morphism $\Spec(R)\to X=\Spec(R^+)$ lifts canonically to a morphism $\Spec(R)\to Z$. Via this morphism, 
the local rings $\cO_{Z,z}$ of $Z$ at $z\in Z$ are regarded as  subrings of the field of fractions of $R$.

\subsubsection{}\label{sY} On the other hand, let $Z'$ be the set of all  pairs $(\mathfrak p, V)$ where $\mathfrak p$ is a prime ideal of $R$ and $V$ is a subring of the residue field $\kappa(\mathfrak p)$ of $\mathfrak p$ satisfying the following conditions (i) -- (iii).

\medskip

(i) $V$ is a valuation ring.

(ii) The image of $R^+$ in $\kappa(\mathfrak p)$ is contained in $V$.

(iii) ${\bar K} \otimes_{\bar A} V=\kappa(\mathfrak p)$.

\begin{sbprop}\label{pV}  We have a bijection between $Z$ and $Z'$ characterized by the following property. If $z\in Z$ corresponds to $(\mathfrak p, V)\in Z'$,
we have $$\cO_{Z, z}= \{f\in R_{\mathfrak p}\;|\; \text{the image of $f$ in $\kappa(\mathfrak p)$ belongs to}\; V\}$$ as a subring of the field of fractions of $R$. 
\end{sbprop}
\noindent This follows from \cite{T} Cor. 3.4.7 taking $\Spec(R)\to \Spec(R^+)$ as $Y\to X$ there.
\\\\
We will  identify the sets $Z$ and $Z'$ via the above canonical bijection.

\subsubsection{}\label{adic1} Via the relation to the theory of adic spaces explained in \ref{adican} below, we use the following notation from adic geometry. 
\\\\For $z=(\mathfrak p, V)\in Z$ and for $f\in R$, let  $|f|(z)$ be the image of $f$ in $\kappa(\mathfrak p)/V^\times=\kappa(\mathfrak p)^\times/V^\times \cup \{0\}$. 
 The set $\kappa(\mathfrak p)/V^\times$ has a natural multiplicative structure and 
is a totally ordered set for the following ordering.  
For $a, b\in \kappa(\mathfrak p)/V^\times$, $a\leq b$ means $V\tilde a\subset V\tilde b$ where $\tilde a$ and $\tilde b$ are representatives of $a$ and $b$ in $\kappa(\mathfrak p)$, respectively. 

\subsubsection{}\label{adican} Let $D^n$ (resp. $\tilde D^n$) be the subset of $Z$ consisting of all points whose images in $\Spec(\bar A)$ are the maximal ideal of $\bar A$ (resp. are not the ideal $(0)$). We endow  $D^n$ and $\tilde D^n$ with the topologies as subspaces of $Z$. 
Then this topology of $D^n$ (resp. $\tilde D^n$) is 
the weakest topology for which the sets
$$\{z\in D^n \;(\text{resp.} \tilde D^n)\;|\; |f|(z) \leq |g|(z)\neq 0\}$$
are open  for all $f,g\in R$ (see \cite{T}, Cor. 3.4.7).\\\\
This description of $\tilde D^n$ as  the subset of $Z'$ with the topology described as above shows that as a topological space, $\tilde D^n$ is identified with the adic space $\Spa(R, R^+)$ if we regard $R$ as a topological ring in which the sets $aR^+$ for non-zero elements $a$ of $\bar A$ form a basis of neighborhoods of $0$ in $R$ (\cite{H1}). 
However in this paper, we do not use the theory of adic spaces. 
\\\\
If $A$ is a complete discrete valuation ring, then $D^n=\tilde D^n$, and this space  coincides, as a topological space,  with the rigid analytic $n$ dimensional unit disc used by Abbes and Saito in \cite{AS1}.  For the theory of rigid analytic spaces, see, for example \cite{A} and \cite{FK}.

\subsubsection{}\label{cODn}

For general valuation rings, we use $D^n$. We endow $D^n$ with the inverse image of the structure sheaf $\cO_Z$ of $Z$ in \ref{sZ}. 

In the theory of adic spaces,  and also in rigid geometry in the case $A$ is a complete discrete valuation ring, the structure sheaf of $\tilde D^n$ is a sheaf of complete rings for the adic topologies. 
But we use  the above structure sheaf of  $D^n$ which  is not  completed. The completion is not a very good operation for general valuation rings.

\subsection{Definition and properties of $X(S, T, I)$.}\label{XSTI} Let $(S,T)$ be a pre-presentation of $B/A$ and let $I$ be a non-zero proper ideal of $\bar A$. 
\\\\
We denote by $D^S$ the locally ringed space defined in the same way as  $D^n$ with $n=\sharp(S)$ replacing the $n$ variables $T_1, \dots, T_n$ by the $n$ variables $y_s$ ($s\in S$). 
 That is, $D^S=D^n$ as soon as we give a  numbering $S=\{s_1, \dots, s_n\}$ of elements of $S$. 
\\\\
We denote by  $X(S, T, I)$ 
 the open subset of $D^S$  consisting of all $z\in D^S$ satisfying the following condition:\\

For each $t\in T$, there is an element $a$ of $I$ such that $|t|(z) \leq |a|_{\bar A}$.

\medskip

\noindent Note that this condition is equivalent to the condition that in the local ring $\cO_{D^S,z}$, $T$ is contained in $I\cO_{D^S,z}$.

\subsubsection{}\label{opencov}

We have an open covering
$$X(S, T, I) =\cup_{a\in I\smallsetminus \{0\}} \;X(S, T, \bar A a).$$

\subsubsection{}\label{I2}
Let $\Phi(L/K)$ be the finite set of all $K$-homomorphisms $L\to \bar K$. Then  
$$\Phi(L/K) \subset X(S, T, I).$$

\noindent The set $\Phi(L/K)$ can be viewed as a subset of $X(S, T, I)$ as follows.  For $\phi\in \Phi(L/K)$, 
$\phi$ is regarded as $(\mathfrak p, V)\in D^S$ (\ref{pV}) where $\mathfrak p$ is the kernel of $\bar K[\{y_s\}_{s\in S}]\to \bar K\;;\; y_s\mapsto \phi(s)$ and $V=\bar A\subset \bar K=\kappa(\mathfrak p)$. 

\section{Upper ramification groups}\label{s:up}
\subsection{Definitions of $\Gal(\bar K/K)^I_{\log}$ and $\Gal(\bar K/K)^I_{\nl}$.}
We define the upper ramification groups by using the notion of ramification bounded by a certain ideal, the latter is defined in \Cref{rambd}.

\begin{sbdefn}\label{defurg}

Let $A$ and $K$ be as in \ref{N1}. For a nonzero proper ideal $I$ of $\bar A$, we define a normal  subgroup $\Gal(\bar K/K)^I_{\log}$ (resp. $\Gal(\bar K/K)^I_{\nl}$) of $\Gal(\bar K/K)$ to be the intersection of the kernels of $\Gal(\bar K/K)\to \Gal(L/K)$ where $L$ ranges over all finite Galois extensions of $K$ in $\bar K$ with ``ramification logarithmically (resp. non-logarithmically) bounded by  $I$". 
\end{sbdefn}

\medskip
\noindent  We will now explain the notion of finite Galois extensions $L/K$ with \tit{ramification $*$-bounded by $I$} for $*=$ logarithmically or non-logarithmically. Let $B$ be the
integral closure of $A$ in $L$.

\begin{sbdefn} For a pre-presentation  $(S, T)$ for $L/K$, we say $(S, T, I)$ \tit{separates $L/K$} 
if for some closed  open subsets $U_{\phi}$ of $X(S, T, I)$ such that $\phi\in U_{\phi}$ for all $\phi\in \Phi(L/K)$, we have

$$X(S, T, I)=\coprod_{\phi \in \Phi(L/K)} U_{\phi}.$$

In other words, $(S, T, I)$ \tit{separates $L/K$} 
 if and only if the embedding $\Phi(L/K) \to X(S, T, I)$ admits a continuous retraction. \end{sbdefn}

\begin{sbprop}\label{S1S2}

Assume that there is an $A$-homomorphism $$A[\{y_s\}_{s\in S_1}]/(T_1)\to A[\{y_s\}_{s\in S_2}]/(T_2)$$ which is compatible with the maps to $B$. Assume further that $(S_1,T_1, I)$ separates $L/K$.\\ Then $(S_2,T_2,I)$ separates $L/K$. 
\end{sbprop}

\begin{pf} Lift this homomorphism to an $A$-homomorphism $A[\{y_s\}_{s\in S_1}] \to A[\{y_s\}_{s\in S_2}]$ which sends $(T_1)$ to $(T_2)$. It induces a morphism of locally ringed spaces $X(S_2,T_2,I)\to X(S_1,T_1, I)$ which is compatible with the maps from $\Phi(L/K)$. This proves \ref{S1S2}.
\end{pf}

\begin{sbcor} (1) Let $(S_1, T_1)$ and $(S_2, T_2)$ be pre-presentations for $L/K$ and assume $S_1\subset S_2$ and $T_1\subset T_2$. If $(S_1, T_1, I)$ separates $L/K$, then $(S_2, T_2, I)$ separates $L/K$.  
\\
(2) If $(S_1, T_1)$ and $(S_2, T_2)$ are presentations for $L/K$, $(S_1, T_1, I)$ separates $L/K$ if and only if  $(S_2, T_2, I)$ separates $L/K$.

\end{sbcor}

\subsection{Ramification bounded by $I$.}
\begin{sbdefn} \label{rambd}
We say the {\it ramification of $L/K$ is non-logarithmically bounded by  $I$} if there is a pre-presentation $(S,  T)$ of $B/A$ such that $(S, T, I)$ separates $L/K$. 
\\\\
We say the {\it ramification of $L/K$ is logarithmically bounded by  $I$} if there is a log smooth extension $A'/A$ of Henselian valuation rings  such that the ramification of $LK'/K'$ is non-logarithmically bounded by $\bar A' I$. 
\\
Here  $K'$ is the field of fractions of $A'$, $L'$ is the composite field $LK'$, and
$\bar A'$ is the integral closure of $A'$ in a separable closure $\bar K'$ of $K'$ which contains $\bar K$ and $L'$ over $K$.

\end{sbdefn}

\begin{sbprop}\label{LK'} Let $*$ be logarithmically or non-logarithmically. Let $L/K$ be a finite Galois extension. Let $A'$ be a Henselian valuation ring which dominates $A$, let $K'$ be the field of fractions of $A'$, let $\bar K'$ be a separable closure of $K'$ which contains $\bar K$, and let $\bar A'$ be the integral closure of $A'$ in $\bar K'$. Assume the ramification of $L/K$ is $*$-bounded by $I$. Then the ramification of $LK'/K'$ is $*$-bounded by $\bar A' I$.

\end{sbprop}

\begin{pf} If $(S, T, I)$ separates $L/K$, $(S, T, \bar A' I)$ separates $LK'/K'$. 
\end{pf}

\begin{sbprop}\label{LL'} Let $*$ be logarithmically or non-logarithmically.
Let $L_1$ and $L_2$ be finite Galois extensions of $K$ in $\bar K$. Assume that the ramification of $L_1/K$ and that of $L_2/K$ are $*$-bounded by $I$. Then for the composite field $L_1L_2$, the ramification of $L_1L_2/K$ is $*$-bounded by $I$. 
\end{sbprop}

\begin{pf}  For $i=1, 2$, let $(S_i, T_i)$ be a pre-presentation for $L_i/K$ such that $(S_i, T_i, I)$ 
separates $L_i/K$. Then   $(S_1 \cup S_2, T_1\cup T_2)$ is a pre-presentation  for $L_1L_2/K$ and $(S_1\cup S_2, T_1\cup T_2, I)$ separates $L_1L_2/K$. This proves the non-log case of  \ref{LL'}. The log case of \ref{LL'} follows from this and from \ref{A'''}. 
\end{pf}

 \begin{para}\label{fin}
 For $*=\log, \nl$, we have defined the upper ramification group $\Gal(\bar K/K)^I_*$.\\
For a finite Galois extension  $L/K$ and for $*=\log, \nl$, we denote the image of $\Gal(\bar K/K)^I_*\to \Gal(L/K)$ by $\Gal(L/K)^I_*$. 

 \end{para}

\begin{sbrem}

 In this paper, we do not write $\Gal(\bar K/K)^I_{\nl}$ simply as $\Gal(\bar
K/K)^I$ (we do not follow \cite{AS1}  concerning this
point). Through the works \cite{V1} and \cite{V2}, we have the
impression that in the study of arbitrary valuation rings, the
logarithmic ramification theory is more natural than the
non-logarithmic theory, and therefore, we have the impression that it is
the logarithmic upper ramification group (not the non-logarithmic one)
that  deserves the simpler notation.

\end{sbrem}

\subsection{Basic Properties.} 
We now give some elementary properties of upper ramification groups.

\begin{sbprop}\label{easy}  We have $\Gal(\bar K/K)^I_{\log}\subset \Gal(\bar K/K)^I_{\nl}$. \\We also have $\Gal(\bar K/K)^I_{\log}\supset \Gal(\bar K/K)^J_{\log}$ and $\Gal(\bar K/K)^I_{\nl}\supset \Gal(\bar K/K)^J_{\nl}$ if $I\supset J$. 

\end{sbprop}

\begin{sbprop}\label{ch0}  Assume that the residue field of $A$ is of characteristic $0$. Then $\Gal(\bar K/K)^I_{\log}=\{1\}$ for every nonzero proper ideal $I$ of $\bar A$. 

\end{sbprop}

\begin{pf} In this case,  any finite Galois extension $L$ of $K$ is a tame extension. Let $A'=B$ be the integral closure of $A$ in $K':=L$. Then $A'/A$ is a log smooth extension and $LK'=K'$. Therefore, the ramification of $L/K$ is logarithmically bounded by $I$.
\end{pf}

\begin{sbprop}\label{LK'2} Let $*=$ log or non-log. Let $A'$, $K'$, $\bar K'$, $\bar A'$ be as in \ref{LK'} and let $I'=\bar A' I$. Then the image of $\Gal(\bar K'/K')_*^{I'}\to \Gal(\bar K/K)$ is contained in $\Gal(\bar K/K)^I_*$. 

\end{sbprop}
\begin{pf} This follows from \ref{LK'}. 
\end{pf}

\medskip 
\noindent In Section \ref{s:pr}, we will prove more properties  of upper ramification groups, including  \Cref{LL'2} and  \Cref{cap} whose proofs are difficult and require additional preparation in sections \ref{s:Sa} and  \ref{s:re}.\\
\begin{sbrem}
We note that the possible usefulness of Zariski-Riemann spaces in the study of upper ramification has appeared in \cite{F}.
\end{sbrem} 

\section{Theory of Saito}\label{s:Sa}
 In  this section, we consider the work \cite{Sa} of T. Saito. 
 Principal ideals $I$ of $\bar A$ appear in these results.  In our ramification theory, we also consider 
 non-principal ideals $I$. We give a complement Proposition \ref{Sapro} that we can deduce  for  non-principal ideals from Saito's results on principal ideals.

\subsection{Connection with \cite{Sa}.}\label{relSa0} We now explain the main point that connects these two theories.
\\\\
Let $L$ be a finite Galois extension of the field of fractions $K$ of $A$ and let $(S, T)$ be a pre-presentation for $L/K$. 
Let  $a$ be a nonzero element of $\m_{\bar A}$ and let  $I$ be the ideal of $\bar A$ generated by $a$.  Let $\mathfrak{Q}^{[D]}= \Spec(\bar A[\{y_s\}_{s\in S}, ta^{-1} (t\in T)])$ and let $\mathfrak{Q}^{(D)}$ be the
normalization of $\mathfrak{Q}^{[D]}$.\\\\ We will consider the scheme
  $$Y(S, T, I):=\bar k \otimes_{\bar A} \mathfrak{Q}^{(D)},$$ where $\bar k$ denotes the residue field of $\bar A$.
  \\
  This scheme is closely related to the ramification theory of Saito in \cite{Sa}.
\\\\  
     We have a surjective continuous closed map
$$X(S, T, I) \to Y(S, T, I)$$
with connected fibers, defined as follows. 
\\ 
Consider the Zariski-Riemann space $Z$ 
 associated with 
 $R=\bar K[\{y_s\}_{s\in S}]$, $R^{+}=\bar A[\{y_s\}_{s\in S}]$ (as in \ref{sZ}). 
 The canonical map $Z\to \Spec(R^{+})$ induces a surjective continuous closed 
map
$$\{z\in Z \;|\; |t|(z)\leq |a|_{\bar A}\;\text{for all}\; t\in T\} \to \mathfrak{Q}^{(D)}.$$
  \medskip \noindent We apply Zariski's main theorem and use the following \cref{zmt} to go to the inverse limit, to show that the fibers of this map are connected. 

The map induced by this map on the inverse images of the closed point of $\Spec(\bar A)$ is $X(S, T, I) \to Y(S, T, I)$.

\begin{sblem}\label{zmt} Let ${(\sX_\la)}_\la$ be a filtered inverse system of topological spaces and let $\sX:= \varprojlim_\la \sX_\la$. Assume that $\sX$ is quasi-compact and that the maps $\pi_\la: \sX \to \sX_\la$ are surjective for all $\la$. Under these assumptions, if $\sX_\la$ is connected for all $\la$, then $\sX$ itself is connected.
\end{sblem}
\begin{pf}
    Suppose that $\sX=U \sqcup V$ is a disjoint union of nonempty open sets $U$ and $V$. Let 
    $U=\cup_{i=1}^n U_{\la_i}$, where each $U_{\la_i}$ is the inverse image of some open set of $\sX_{\la_i}$.

    There exists $\la_0$ (we will call it $0$ for simplicity) such that $U=\pi_0^{-1}(U_0)$ and $V=\pi_0^{-1}(V_0)$, where $U_0, V_0$ are open sets of $\sX_0$. Note that they are both nonempty since $U, V$ are nonempty.
    By $U \cup V =\sX$ (resp. $U \cap V = \phi$), and by the surjectivity of $\pi_0$, we have $U_0 \cup V_0 =\sX_0$ (resp. $U_0 \cap V_0 = \phi$).
    
\end{pf}

\subsection{Complete intersection.}
\label{Sa4}
   We say $(S, T)$ is of complete intersection if $A\to B':=A[\{y_s\}_{s\in S}]/(T)$ is locally of complete intersection. 
  Note that the map $A\to B'$ is flat because it is injective (the composition $A\to B'\to B$ is injective)   and any torsion free module over a valuation ring is flat.  
\\\\
In the rest of the section 5, assume that $(S, T)$ is of complete intersection and assume that $K[\{y_s\}_{s\in S}]/(T) \overset{\cong}\to L$.

\noindent We provide some explanation and further context for \cite{Sa} below.\\\\
Let $\cR$ be a  subring of $\bar{A}$ satisfying the following conditions (i) -- (iii):\\
(i) $\cR$ is normal and finitely generated over $\Z$.\\
(ii) The element $a \in \bar{A}$ generating the ideal $I$ is in $\cR$.\\
(iii) $\cR[\{y_s\}_{s\in S}]$ contains the image of $T \subset A[\{y_s\}_{s\in S}]$ in $\bar{A}[\{y_s\}_{s\in S}]$. \\
We regard $T$ as a subset of $\cR[\{y_s\}_{s\in S}]$.
\\\\
We apply \cite[section 2.1 ]{Sa} by taking $$X:= \Spec(\cR), Y:= \Spec(\cR[\{y_s\}_{s\in S}]/(T)), Q:= \Spec(\cR[\{y_s\}_{s\in S}]).$$ 
\\
Then the schemes $Q^{[D]}$ and $Q^{(D)}$ in  \cite[2.1]{Sa} for the divisor $D$ on $X$ defined by $a$ are as follows.  
$$Q^{[D]}:=\Spec(\cR[\{y_s\}_{s\in S}, ta^{-1} (t \in T)])$$
and $Q^{(D)}$ is the normalization of $Q^{[D]}$.
\\The scheme $\mathfrak{Q}^{(D)}$ in \ref{relSa0} is the inverse limit of $Q^{(D)}$ as $\cR$ varies, and hence, $Y(S,T,I)$ from \ref{relSa0} is the inverse limit of $\bar{k} \otimes_{\cR} Q^{(D)}$.
\\\\
The following lemma is a consequence of 
\cite{STACKS}.
\begin{sblem}\label{bigR} For a sufficiently large $\cR$, the morphism $Y \to X$ is flat and locally of complete intersection.
\end{sblem}

\noindent This allows us to obtain the following results
for sufficiently large $\cR$, as corollaries of  \cite[Proposition 3.1.2 (1) and Lemma 3.2.4]{Sa}, respectively.
\begin{sbcor}\label{cor1s}
 The canonical map $\Phi(L/K) \to \pi_0(\bar{k} \otimes_{\cR} Q^{(D)})$ is surjective.    
\end{sbcor} 
\begin{sbcor}\label{cor2s}
    If $\cR'$ is a normal subring of $\bar{A}$ finitely generated over $\cR$, and if $Q'^{(D)}$ is the analog of $Q^{(D)}$ for $\cR'$, the map 
$\pi_0(\bar{k} \otimes_{\cR'} Q'^{(D)}) \to \pi_0(\bar{k} \otimes_{\cR} Q^{(D)})$ is bijective.
\end{sbcor}

\subsection{Principal ideals.}

\subsubsection{}\label{Sasep} From  \cref{cor1s} and \cref{cor2s}, we see that $Y(S, T, I)$ is a finite disjoint union of certain connected closed open subspaces (connected components) of $Y(S, T, I)$. The map $\Phi(L/K)\to \pi_0(Y(S, T, I))$ is surjective. 
\\
 Consequently, the following three conditions are equivalent:\\

(i) The number of connected components of $Y(S, T, I)$ is $[L:K]$.

(ii) $Y(S, T, I)$ is a disjoint union of  closed open  subspaces $U_{\phi}$ ($\phi\in \Phi(L/K)$) such that for each $\phi\in \Phi(L/K)$, the image of $\phi$ in $Y(S, T, I)$ is contained in $U_{\phi}$ .

(iii) The map $\Phi(L/K) \to \pi_0(Y(S, T, I))$ is bijective. 
\\\\
 The ramification of $Y \to X$ at $\bar{A}$ is bounded by $I=a\bar{A}$ in the sense of \cite[3.2.3]{Sa} if and only if these equivalent conditions are satisfied.  
 \\\\
By  \cite[Theorem 3.2.6]{Sa}, we have: 
\begin{sbprop}\label{SsTt2}
There is a principal nonzero  ideal $I_1$ of $\bar A$ such that for any principal nonzero proper  ideal $I$ of $\bar A$,   the ramification of $Y \to X$ at $\bar{A}$ is bounded by $I$ if and only if $I\subsetneq I_1$. 
\end{sbprop}

\subsubsection{}\label{STS'T'} By  \cite[Proposition 3.1.2 (2) and Lemma 3.2.4]{Sa} we have the following:\\\\
Let $I$ be a principal nonzero proper ideal of $\bar A$. 
Let  $L'$ be a finite Galois extension of $K$ contained in $L$, let $S'\subset S$, $T'\subset T$, and assume  $S'\subset L'$ and $K[\{y_s\}_{s\in S'}]/(T') \overset{\cong}\to L'$ and that $(S',T')$ is of complete intersection. Then $\pi_0(Y(S', T', I))$ is the pushout of $$\Phi(L'/K) \leftarrow \Phi(L/K) \to \pi_0(Y(S, T, I))$$ in the category of sets. 
\\\\
We now consider our space $X(S,T,I)$. Recall that (\ref{relSa0}) we have a surjective continuous closed map $X(S, T, I) \to Y(S, T, I)$ with connected fibers.
\begin{sblem}\label{123}
Let $f: \sX\to \sY$ be a surjective continuous closed map of topological spaces whose fibers are connected. Then we have a bijection $C\mapsto f^{-1}(C)$  from the set of all closed open subsets of $\sY$ onto the set of all closed open subsets of $\sX$. 
\end{sblem}
\begin{pf} Assume that $\sX$ is the disjoint union of two closed  open subsets $C_1$ and $C_2$. For $i=1, 2$, the image $C'_i$ of $C_i$ in $\sY$ is closed. We prove that $C'_1\cap C'_2=\emptyset$ and that $C_i=f^{-1}(C'_i)$. Let $y\in C'_1\cap C'_2$. Let $F\subset \sX$ be the fiber of $y$. Then  $F$ is the disjoint union of its closed open subsets $F\cap C_i$. Since $F$ is connected, we have $F\cap C_i=\emptyset$ for some $i$. But this contradicts the fact that $y$ is in the image of $C_i$. Therefore, we have $C'_1\cap C'_2=\emptyset$. It follows that $C_i=f^{-1}(C'_i)$. \end{pf}

\subsubsection{}
Applying this to $\sX=X(S, T, I)$ and $\sY=Y(S, T, I)$, we see
 that for $I= a{\bar A}$, $(S, T, I)$  separates $L/K$ in our sense if and only if the ramification of $\cR[\{y_s\}_{s\in S}]/(T)$  over $\cR$ at $\bar A$ is bounded by $I$ in the sense 
 of  Saito.
 \\ 
We also have 
$$\pi_0(X(S, T, I))\overset{\cong}\to \pi_0(Y(S, T, I)).$$
\\

\noindent Now we shift our focus to include non-principal ideals.

\subsection{Non-principal ideals.}

\noindent Discussions in \ref{prep1} and \ref{prep2} lead us to our general result \Cref{Sapro} below.

\begin{sbprop}\label{Sapro} Let $(S, T)$ be a pre-presentation for $L/K$ which is of complete intersection such that $K[\{y_s\}_{s\in S}]/(T) \overset{\cong}\to L$ and let $I$ be a nonzero proper ideal of $\bar A$. 
\begin{enumerate}[(1)]
\item The space $X(S, T, I)$ is a finite disjoint union of its connected components. The map $\Phi(L/K)\to \pi_0(X(S, T, I))$ is surjective. 

\item The tuple $(S, T, I)$ separates $L/K$ if and only if $(S, T, J)$ separates $L/K$ for all principal nonzero subideals $J$ of $I$. 

\item  Replace $Y(S, T, I)$ with $X(S, T, I)$ in \ref{Sasep}. Then the statements in \ref{Sasep} remain true.

\item Let $I_1$ be as in \cref{SsTt2}. Then for any nonzero proper ideal $I$ of $\bar A$, $(S, T, I)$ separates $L/K$ if and only if $I\subsetneq I_1$. 

\item 
Let  $L'$ be a finite Galois extension of $K$ contained in $L$, let $S'\subset S$, $T'\subset T$, and assume  $(S',T')$ is of complete intersection and that 
$S'\subset L'$ and $K[\{y_s\}_{s\in S'}]/(T') \overset{\cong}\to L'$. Then 
$\pi_0(X(S', T', I))$ is the pushout of $$\Phi(L'/K) \leftarrow \Phi(L/K) \to \pi_0(X(S,T,I))$$ in the category of sets. 
\end{enumerate}

\end{sbprop}

\subsubsection{$\text{Map}_c(\sZ, \F_2)$.}\label{prep1}
For a topological space $\sZ$, let $\text{Map}_c(\sZ, \F_2)$ be the set of all continuous maps from $\sZ$ to the discrete field $\F_2=\Z/2\Z$. Then $\text{Map}_c(\sZ, \F_2)$ is identified with the set of all closed open subsets of $\sZ$ (an element of $\text{Map}_c(\sZ, \F_2)$ gives such a subset of $\sZ$ as the inverse image of $\{0\}$).   $\text{Map}_c(\sZ, \F_2)$ is a finite set if and only if 
  $\sZ$ is a finite disjoint union of its connected components.
If $\sZ$ is such a space,   via the canonical map $\sZ \to \pi_0(\sZ)$, we have a bijection $\text{Map}(\pi_0(\sZ), \F_2) \to \text{Map}_c(\sZ, \F_2)$ and $\pi_0(\sZ)$ is identified with the set of all ring homomorphisms $\text{Map}_c(\sZ, \F_2)\to \F_2$. 
\subsubsection{Proof of \Cref{Sapro}.}\label{prep2}
Assume $(S, T)$ is of complete intersection and assume $K[\{y_s\}_{s\in S}]/(T) \overset{\cong}\to L$. Let $I$ be a nonzero proper ideal of $\bar A$. Recall that $X(S, T, I)$ is the union of the open sets $X(S, T, J)$ where $J$ ranges over all principal nonzero subideals of $I$ (\ref{opencov}).\\ Hence, we have 
 $\text{Map}_c(X(S,T, I), \F_2) = \varprojlim_J \text{Map}_c(X(S,T,J), \F_2)$ where $J$ ranges over all principal  subideals of $I$. Since the canonical map $\Phi(L/K)\to \pi_0(X(S, T, J))$ is surjective, as a quotient of $\Phi(L/K)$, $\pi_0(X(S, T, J))$ is independent of a sufficiently large principal  subideal $J$ of $I$. 
 
 Therefore,  
the map $\text{Map}_c(X(S,T, I), \F_2) \to \text{Map}_c(X(S,T, J), \F_2)$ is an isomorphism for such $J$. In particular, $\text{Map}_c(X(S, T, I), \F_2)$ is finite. Consequently, we have (1) and (2) in  \Cref{Sapro}, and we can obtain (3), (4), (5) in it from the results \ref{Sasep},  \ref{SsTt2}, \ref{STS'T'}, respectively.

\section{Applications of the works  \cite{V1}, \cite{V2}, \cite{V3}}\label{s:re}
\cite{V1, V2, V3} deal with degree $p$ extensions in positive residue characteristic $p$.
 \cite{V1} treats the equal characteristic case when $L/K$ is either defectless with valuation of any rank or has defect and rank $1$ valuation. \cite{V2} treats the mixed characteristic case without any restriction on the rank of the valuation. \cite{V3} completes the study of equal characteristic case by removing aforementioned restrictions  on rank, we also provide some additional proofs in this work.
\\\\ 
 In this section, we prove the following results as applications of \cite{V1, V2, V3}. 
\begin{thm}\label{prime} Let $L/K$ be a cyclic extension of degree a prime number. Let $S$ be a finite subset of $B$. Then there is $s\in B$ such that $S\subset A[s]$. 

\end{thm}
\noindent (Proof in \ref{pf*}-\ref{elem}).
\begin{thm}\label{thm2}  Let $L$ be a finite Galois extension of $K$ and let $B$ be the integral closure of $A$ in $L$. Then for any finite subset $S$ of $B$, there is an
$A$-subalgebra $B'$ of $B$ which is  free of finite rank as an $A$-module  and of complete intersection over $A$ such that $S\subset B'$. 

\end{thm}
\noindent (Proof in \ref{pfthm2}.)
\begin{cor}\label{thm2cor} Let the notation be as in  \Cref{thm2}.  If $B$ is a finitely generated $A$-module, $B$ is of complete intersection over $A$. 

\end{cor} 

\begin{thm}\label{appl} Let $L/K$ be a finite Galois extension. Then there is a log smooth extension of Henselian valuation rings $A'/A$ with field of fractions $K'$ such that  $e(LK'/K')=1$. Here $e$ denotes the ramification index (\ref{N1}). We can take such $A'$ such that $A=A_0\subset A_1\subset \dots \subset A_n=A'$ where $A_i/A_{i-1}$ for $i=1$ (resp. $2\leq i\leq n$) is a log smooth extension of Henselian valuation rings of type 1 (resp. 2) (we do not need a log smooth extension of type 3) (\ref{types}). 
 \end{thm}
 \noindent (Proof in \ref{pfappl}.)
 
\begin{thm}\label{dlequiv}  Let $L/K$ be a finite Galois extension. Then the  following three conditions are equivalent.
\\\\
(i) $L/K$ is defectless. 
\\
(ii) There is a log smooth extension $A'/A$ of Henselian valuation rings with field of fractions $K'$ such that the integral closure $B'$ of $A'$ in $LK'$ is a finitely generated $A'$-module. 
\\
(iii) There is a log smooth extension $A'/A$ of Henselian valuation rings with field of fractions $K'$ such that the integral closure $B'$ of $A'$ in $LK'$ is finite flat and complete intersection over $A'$.
\\\\
We still have the equivalence of conditions when we replace the log smoothness of $A'$ over $A$ in (ii) and (iii) by the stronger property that there is a sequence $A=A_0\subset A_1\subset \dots \subset A_n=A'$ such that $A_i/A_{i-1}$ for $i=1$ (resp. $2\leq i\leq n$) is a log smooth extension of Henselian valuation rings of type 1 (resp. 2). 
\end{thm}
\noindent (Proof in \ref{7.5pf}.)\\\\
 In \ref{rev1}-\ref{zetap}, we review some material from \cite{V1, V2, V3}.
In \ref{rev1}--\ref{pf*}, we assume that the residue field of $A$ is of positive characteristic $p$, and we denote by $L$ a cyclic extension of $K$ of degree $p$. 
Let $B$ be the integral closure of $A$ in $L$ and let $l$ be the residue field of $B$. 
\subsection{Review of special ideals.}
\subsubsection{Ideal $\sH$.}\label{rev1} 
 The important ideal $\sH$ of $A$ introduced in \cite{V1} and \cite{V2} associated to $L/K$  is given as follows.
\\
 Let $G=\Gal(L/K)$. For $\sigma \in G\smallsetminus \{1\}$, let $\js$ be the ideal of $B$ generated by all elements of the form $\sigma(x)/x-1$  ($x\in L^\times$). This ideal $\js$ of $B$ does not depend on $\sigma$. 
\\
The ideal $\sH$ of $A$ is the ideal generated by  $\{N_{L/K}(b)\;|\;b\in \js\}$, where $N_{L/K}$ is the norm map $L\to K$.

\subsubsection{Alternate description of $\sH$.}\label{rev2}  

\noindent Let $\chi\in H^1(K, \Z/p\Z)=\Hom_{\text{cont}}(\Gal(\bar K/K), \Z/p\Z)$ be a nonzero element which gives the cyclic extension $L/K$. The ideal $\sH$ is also described as follows. 
\\\\
Assume $K$ is of characteristic $p>0$. Then by Artin-Schreier theory we have an isomorphism $ H^1(K, \Z/p\Z)\cong K/\{x^p-x\; |\; x \in K\}$  and $\sH$ is the ideal of $A$ generated by $1/f$ for all nonzero elements $f\in K$ which give $\chi$ via the above isomorphism. (See \cite[Theorem 0.3]{V1} and \cite[\S4]{V3}.)
\\\\
Assume $K$ is of mixed characteristic $(0, p)$ and assume that $K$ contains a primitive $p$-th root $\zeta_p$ of $1$. Then we have an isomorphism $H^1(K, \Z/p\Z)\cong K^\times/(K^\times)^p$ by Kummer theory. If there is an element $a$ of $1+\m_A$ which gives $\chi$ via this isomorphism, $\sH$ is the ideal of $A$ generated by $(\zeta_p-1)^p(a-1)^{-1}$ for all such elements $a$. If there is no element of $1+\m_A$ which gives $\chi$ via this isomorphism, then $\sH=(\zeta_p-1)^pA$. (See \cite[Theorem 1.3]{V2}.)
\\\\
If $K$ is of mixed characteristic $(0, p)$ and does not contain $\zeta_p$, $\sH= A\cap \sH'$ where $A'$ is the integral closure of $A$ in $K'=K(\zeta_p)$ 
and $\sH'$ is the ideal of $A'$ associated to $LK'/K'$. (See \cite[Theorem 7.8]{V2}.)

\subsubsection{$\sH$, $\js$, and defect.}\label{rev3}  

Let the notation be as in \ref{rev1}. The following conditions (i) -- (iii) are equivalent.
\\
(i)  $L/K$ is defectless.
\\
(ii) The ideal $\sH$ of $A$ is principal.
\\
(iii) 
The ideal $\js$ of $B$ is  principal. 
\\
(See \cite[Corollary 4.5]{V1} , \cite[Corollary 4.4]{V2} , \cite[\S3 and \S4]{V3}.)

\subsection{Further review of the defectless case.}

\noindent In this subsection, we describe in further detail the case when $L/K$ is defectless.

\subsubsection{Refined Swan conductor}\label{dl0} 

(See \cite[Theorem 0.5]{V1} , \cite[Theorem 1.5]{V2}.)
Assume that $L/K$ is defectless and $L/K$ is not unramified. Then there is a canonical nonzero $A$-homomorphism $$\text{rsw}(\chi): \sH\to k \otimes_A \Omega^1_A(\log)$$
(called the refined Swan conductor) characterized by the following property. Let $\sigma$ be an element of $\Gal(L/K)$ such that $\chi(\sigma)=1\in \Z/p\Z$. Then for $x\in L^\times$, it sends $N_{L/K}(\sigma(x)x^{-1}-1)$ to the class of $d\log(N_{L/K}(x))$.

In \cite{V1} and \cite{V2}, 
  an $A$-homomorphism from $\sH$ to a certain bigger quotient of $\Omega^1_A(\log)$ is defined without assuming $L/K$ is defectless, but in this paper, we use this map $\text{rsw}(\chi)$ induced by it under the defectless situation. Since $\sH$ is principal in the defectless case, this homomorphism $\text{rsw}(\chi)$ is regarded as an element of $k \otimes_A \Omega^1_A(\log)\otimes_A \sH^{-1}$ where $\sH^{-1}$ is the inverse fractional ideal of $\sH$.

\subsubsection{Best $f$ \cite[\S3.2]{V3}.}\label{dl1}  Assume $K$ is of positive characteristic $p$ and assume $L/K$ is defectless. Then $L=K(\alpha)$ with $\alpha^p-\alpha=f$ for $f\in K^\times$ satisfying one of the following conditions (i) -- (iii).
\\\\
(i) $f\notin A$ and the class of $f$ in $K^\times/A^\times$ is not a $p$-th power.
\\
(ii) $f\notin A$ and $f=ug^{-p}$ for some $g\in K^\times$ and $u\in A$ such that the residue class of $u$ is not a $p$-th power. 
\\
(iii) $f\in A$ and the residue class of $f$ does not belong to $\{x^p-x\;|\; x\in k\}$. 
\\\\
In the case (i), the ramification index of $L/K$ is $p$. In the case (ii), the residue class of $B$ is a purely inseparable extension of $k$ generated by the $p$-th root of the residue class of $u$.  In the case (iii), $L/K$ is unramified. In any of these cases (i)--(iii), $\sH$ is generated by $f^{-1}$.  In the cases (i) and (ii), $\rsw(\chi)$ sends $f^{-1}$ to $d\log(f)$. This follows from $\sigma(\al)\al^{-1}-1=(\al+1-\al)\al^{-1}=\al^{-1}$ and  $N_{L/K}(\al)=f$.

\begin{sbrem}\label{7794} 
   Recall the exact sequence  
    $\begin{matrix}  0 & \to &  \Omega^1_k & \to & k \otimes_A \Omega^1_A(\log) & \to & k \otimes_\Z \Gamma_A & \to 0\end{matrix}$ from  \ref{ses}.

    \noindent In the case (i), the class of $f$ in $\Gamma_A$ is non-trivial. Since the image of $\dl f$ in the above sequence is non-trivial, we see that $\rsw$ is a non-zero map.

    \noindent In the case (ii), $\dl f=\dl u$ is non-trivial since the residue class of $u$ is not a $p$-th power. Once again, $\rsw$ is a non-zero map.
\end{sbrem}

\subsubsection{Best $h$ \cite[\S3.3]{V3}.}\label{dl2}
Assume $A$ is of mixed characteristic $(0, p)$ and assume $A$ contains a primitive $p$-th root of $1$. Assume $L/K$ is  defectless.  Then $L=K(h^{1/p})$ with $a\in K^\times$ satisfying one of the following  conditions (i) -- (v).
\\\\
(i) The class of $h$ in $K^\times/A^\times$ is not a $p$-th power.
\\
(ii) $h\in A^\times$ and the residue class of $h$ is not a $p$-th power.
\\
(iii) $h= 1+ b$ for some $b\in K^\times$ such that $b, (\zeta_p-1)^pb^{-1}\in \m_A$ and such that the class of $b$ in $K^\times/A^\times$ is not a $p$-th power. 
\\
(iv) $h=1+b$ for some $b\in K^\times$ such that $b, (\zeta_p-1)^pb^{-1}\in \m_A$ and such that $b=g^pu$ for some $g\in K^\times$ and $u\in A$ such that the residue class of $u$ is not a $p$-th power. 
\\
(v) $h=1+(\zeta_p-1)^pb$ where $b\in A$ and the residue class of $b$ does not belong to $\{x^p-x\;|\; x\in k\}$. 
\\\\
In the cases (i) and (iii), the ramification index of $L/K$ is $p$. In the case (ii) (resp. (iv)), the residue class of $B$ is a purely inseparable extension of $k$ generated by the $p$-th root of the residue class of $h$ (resp. $u$). In the case (v), $L/K$ is unramified. In the cases (i) and (ii), $\sH=A(\zeta_p-1)^p$. In the cases (iii) and (iv), $\sH=A(\zeta_p-1)^pb^{-1}$. In the case (v), $\sH=A$. In the cases (i) and (ii), $\text{rsw}(\chi)$ sends $(\zeta_p-1)^p$ to $d\log(h)$. In the cases (iii) and (iv), $\text{rsw}(\chi)$ sends $(\zeta_p-1)^pb^{-1}$ to $d\log(b)$. This follows from \cite{V2} Theorem 1.5(i).

\begin{sbrem}Analogous to \cref{7794}, $\rsw$ is a non-zero in the cases (i) -- (iv) when we have a Kummer extension $L/K$.
The non-Kummer case when $L/K$ is not unramified follows from this.
\end{sbrem}

\subsection{Further review of the defect case.}
\subsubsection{Filtered union.}\label{filuni} In \cite{V1} and \cite{V2}, 
Theorem \ref{prime} for a  defect extension  $L/K$ is proved in the following cases:\\
\\
(i) \cite[Theorem 5.1]{V1}: $K$ is of characteristic $p$. 
\\
(ii)  \cite[Theorem 5.1]{V2}:  $A$ is of mixed characteristic and $A$ contains a primitive $p$-th root $\zeta_p$ of $1$.

\subsubsection{Review of the mixed characteristic case.}\label{zetap} We review \cite[Theorem 5.1, Lemma 5.3, and the proof of Lemma 6.10]{V2} for a defect extension $L/K$. Assume $A$ is of mixed characteristic $(0,p)$ and assume $\zeta_p\in A$. Let $$\sS=\{\alpha\in L^\times\;|\;\alpha^p \in 1+m_A, L=K(\alpha)\}.$$ Then $\sS$ is not empty and we have:
\\\\
(1) $\js$ is generated by elements of the form $(\zeta_p-1)(\alpha-1)^{-1}$ where $\alpha$ ranges over all elements of $\sS$. 
\\
(2) $B=\cup_{\alpha\in\sS} A[\alpha']$ where $\alpha'$ denotes an element of $B^\times \cap A(\alpha-1)(\zeta_p-1)^{-1}$.\\ (Then $A[\alpha']$ depends only on $\alpha$ and it does not depend on the choice of $\alpha'$). 
\\
(3) For $\alpha_1,\alpha_2\in \sS$, we have  $A[\alpha'_1]\subset A[\alpha'_2]$  if $B(\alpha_2-1) \subset B(\alpha_1-1)$. 
\\\\
These (2) and (3) prove the case (ii) in \ref{filuni}.

\subsection{Preparation for \Cref{thm1} and  \Cref{appl}. }
 \Cref{HH'} will be used in the proof of the theorems \ref{thm1} and  \ref{appl}; here $L/K$ may or may not have non-trivial defect.

\begin{sbprop}\label{HH'}  Let $A'/A$ be a log smooth extension of Henselian valuation rings, let $K'$ be the field of fractions of $A'$, and  let $L'=LK'$. Let $\sH'$ be the ideal of $A'$ associated to  $L'/K'$. Then  $\sH'=A'\sH$. 
\end{sbprop}

\begin{pf}  Assume first that $e(L/K)=1$.  

In this case, the ideals  $\is:=(\{ (\si-1)(b) \mid b \in B \})$ and $\js$ for $L/K$ are equal (note that both ideals are independent of the choice of a generator $\si$ of the Galois group of $L/K$). Basic properties of the ideals $\js, \is$ of $B$ are explored in \cite{V1, V2}.

By  \Cref{124},  we also have $e(L'/K')=1$. So, the ideals $\js$ and $\is$ for $L'/K'$ are equal as well.
Again by \Cref{124}, $B'$ is generated by $B$ as an $A'$-algebra. That is, any $b' \in B'$ can be written as a finite sum $\sum ba'$, where $b \in B$ and $a' \in A' \subset K'$. 
For any 
generator $\sigma$ of $\Gal(L'/K')$,  $(\si-1)(\sum ba')=\sum (\si-1)(ba')=\sum a' (\si-1)(b)$.  

Thus, the ideal $\is$ associated with $L'/K'$ is generated by the ideal $\is$ 
associated with $L/K$. 
Consequently, the ideal $\js$ associated with $L'/K'$ is generated by the ideal $\js$ associated with $L/K$. Since $\sH$ is generated by norms of elements of $\js$, this proves $\sH'=A'\sH$. 
\\\\
Assume next that $e(L/K)>1$. 
Clearly,  the extension $L/K$ is defectless. Since $K'/K$ is log smooth, $LK'/L$ is also log smooth. Log smooth extensions are defectless. Consequently, $LK'/K$ is defectless. And hence, $LK'/K'$ is also 
defectless. 
\\
Since, $e(L/K)>1$ and we are in the degree $p$ case, by \Cref{tame'}, we note that $LK'/K'$ is not unramified.
\\\\
We clearly have $\sH'\supset A'\sH$,  and by the construction of the refined Swan conductor as described in \cref{dl0}, we have a commutative diagram

$$\begin{matrix}  \sH & \to & k \otimes_A \Omega^1_A(\log)\\
\downarrow&&\downarrow\\ \sH' & \to & k' \otimes_{A'} \Omega^1_{A'}(\log)
\end{matrix}$$
where $k'$ is the residue field of $A'$ and the upper (resp. lower) horizontal arrow is the refined Swan conductor of $L/K$ (resp. $LK'/K'$). 
\\\\
Since $k \otimes_A \Omega^1_A(\log)\to  k' \otimes_{A'} \Omega^1_{A'}(\log)
$ is injective, the composition 
$\sH\to k \otimes_A \Omega^1_A(\log)\to  k' \otimes_{A'} \Omega^1_{A'}(\log)$ is not zero. But if $\sH'\neq A'\sH$, that is, if $A'\sH \subset \m_{A'}\sH'$, then the composition $\sH \to \sH' \to k' \otimes_{A'} \Omega^1_{A'}(\log)$ must be zero. This proves $\sH'=A'\sH$.
\end{pf}

\noindent  \Cref{dl3} below will be used in the proof of \Cref{appl}. 
\begin{sblem}\label{dl3} 
Assume $L/K$ is defectless and is not unramified. 
\\
Then   $e(L/K)=p$ (resp. $[l:k]=p$) if and only if the image of $\text{rsw}(\chi)$ under $k \otimes_A \Omega^1_A(\log) \to k \otimes_{\Z} \Gamma$ is non-trivial (resp. trivial). 
\end{sblem}
\begin{pf}This follows from \ref{dl1} and  \ref{dl2} (in the mixed characteristic case, we are reduced to the case 
$\zeta_p\in A$ by \ref{HH'}). 
\end{pf}

\subsection{Proof of \Cref{prime} - Part I}

\subsubsection{}\label{filuni2}

In \ref{lemuni}--\ref{pf*},  we prove 
\\\\
(*)  \Cref{prime} for a defect extension $L/K$ in the mixed characteristic case  without assuming $\zeta_p\in A$, 
by reducing it to the case $\zeta_p\in A$ treated in \ref{zetap}.

\begin{sblem}\label{lemuni}
Let $y\in B$, $z\in A[y]$, and assume $\sigma(z)-z\in (\sigma(y)-y)A[y]^\times$ for any non-trivial element $\sigma$ of $\Gal(L/K)$. Then $A[y]=A[z]$.

\end{sblem}

\begin{pf}

We may assume $y\notin A$.
\\
For a free $A$-module $E$ in $L$ such that $K\otimes_A E=L$, let $E^*=\{x\in L\;|\; Tr_{L/K}(xE)\subset A\}$, where $Tr_{L/K}$ is the trace map $L\to K$. We have $(E^*)^*=E$. 

Let $f$ (resp. $g$) be the monic irreducible polynomial of $y$ (resp. $z$) over $K$ and let $f'$ (resp. $g'$) be its derivative. Then $A[y]^*=f'(y)^{-1}A[y]$ and $A[z]^*=g'(z)^{-1}A[z]$. 

We have $f'(y)=\prod_{\sigma} (y-\sigma(y))$, $g'(z)= \prod_{\sigma} (z-\sigma(z))$ where $\sigma$ ranges over all non-trivial elements of $\Gal(L/K)$. Hence, $A[z]^*=g'(z)^{-1}A[z]\subset f'(y)^{-1}A[y]=A[y]^*$. It follows that $A[z]=(A[z]^*)^*\supset (A[y]^*)^*=A[y]$.

\end{pf}

\subsubsection{}\label{pf*} Now we prove (*) from \ref{filuni2}. 
\\\\
Let $K_1=K(\zeta_p)$, $L_1=L(\zeta_p)$, let $A_1$ be the integral closure of $A$ in $K_1$ and let $B_1$ be the integral closure of $A_1$ in $L_1$. 
\\\\
We have $\Gal(L_1/K_1) \overset{\cong}\to \Gal(L/K)$, $\Gal(L_1/L)\overset{\cong}\to \Gal(K_1/K)$. We will regard these isomorphisms as identifications. In particular, denote the generator of $\Gal(L_1/K_1)$ corresponding to $\sigma\in \Gal(L/K)$ by the 
same notation $\sigma$.

Let  $\sS\subset L_1$ be as in \ref{zetap} for $L_1/K_1$. Let $\sT$ be the set of all $\alpha\in \sS$ such that the ideal 
$B_1(\zeta_p-1)(\alpha-1)^{-1}$ of $B_1$ is generated by an element of $B$ (that is, by an element of $A$; note that $\Gamma_A=\Gamma_B$). 
\\\\
For $\alpha\in \sS$, let $\alpha'\in L_1$ be as in \ref{zetap} (using the present $L_1/K_1$ as $L/K$ in  \ref{zetap}). 
\\
We prove (*) in \ref{filuni2} assuming the claims 3 and 4 stated below (these will be proved later in \cref{pfcl3} and \cref{pfcl4}, respectively). 
\\
{\bf Claim 3.} Let $\gamma$ be a nonzero element of $\js$ such that $\gamma$ divides $\zeta_p-1$. Then there is an element of $\alpha\in \sT$ such that $B_1\gamma = B_1(\zeta_p-1)(\alpha-1)^{-1}$.
\\
{\bf Claim 4.}  Assume $\alpha\in \sT$. Then $A_1[\alpha']= A_1[\alpha'']$ for some $\alpha''\in B$. 

Note that $J_{1,\sigma}$ (defined as the $L_1/K_1$-version of $\js$) is equal to $B_1\js$ (\cite[Proposition 7.7]{V2} ).
Hence by claims 3, 4, and \cref{zetap} applied to the extension  $L_1/K_1$,  we have $B_1= \cup_{\alpha\in \sT} A_1[\alpha'']$. By taking the $\Gal(L_1/L)$-fixed parts, we have $B= \cup_{\alpha\in \sT} A[\alpha'']$.

Furthermore, for $\alpha_1, \alpha_2\in \sT$, if $B_1(\alpha_2-1)\subset B_1(\alpha_1-1)$, then we have $A_1[\alpha''_1]\subset A_1[\alpha''_2]$ by \cref{zetap}. By taking the $\Gal(L_1/L)$-fixed parts, we have  $A[\alpha''_1]\subset A[\alpha''_2]$.
 These prove (*).

\subsubsection{Proof of Claim 3.}\label{pfcl3}
\begin{proof}
 Take $\alpha\in \sS$ such that $B_1\gamma \subset B_1(\zeta_p-1)(\alpha-1)^{-1}$. If this inclusion is an equality, then $\alpha\in \sT$. Assume that this is a strict inclusion. Then take any $\delta\in \m_{A_1}$  such that $B_1\delta= B_1(\zeta_p-1)\gamma^{-1}$ (such $\delta$ exists because $\Gamma_{A_1}=\Gamma_{B_1}$). Then $\alpha(1+\delta)$ 
belongs to $\sS$ and $B_1(\alpha(1+\delta)-1)= B_1\delta$ and hence, $B_1(\zeta_p-1)(\alpha(1+\delta)-1)^{-1}=B_1\gamma$. Thus $\alpha(1+\delta)\in \sT$ and this element has the property of $\alpha$ in claim 3. 
\end{proof}

\subsubsection{Proof of Claim 4.}\label{pfcl4}

Let  $\kappa: \Gal(L_1/L)=\Gal(K_1/K)\to (\Z/p\Z)^\times$ be the homomorphism $\tau \mapsto r$, $\tau(\zeta_p)= \zeta_p^r$. 
We have \begin{align*}
H^1(K, \Z/p\Z) & \overset{\cong}\to H^1(K_1, \Z/p\Z)^{\Gal(K_1/K)}\\
& = \{\chi\in H^1(K_1, \Z/p\Z(1))\;|\; \tau(\chi)= \kappa(\tau)\chi\;\forall\;\tau\in \Gal(K_1/K)\} \\
& = \{a\in (K_1)^\times/((K_1)^\times)^p\;|\; \tau(a)= a^{\kappa(\tau)}\;\forall\;\tau\in \Gal(K_1/K)\}. \end{align*}
Let $\tilde \kappa(\tau)\in \Z$ be a lifting of $\kappa(\tau)\in (\Z/p\Z)^\times$.
Then $L_1= K_1(\alpha)$, $\alpha^p=a\in (K_1)^\times$, $\tau(\alpha)= \alpha^{\tilde \kappa(\tau)}c_{\tau}$ for $\tau\in \Gal(K_1/K)$ and for $c_{\tau}\in (K_1)^\times$. 
Because $L_1/K_1$ is a defect extension,  
we can take $a\in 1+\m_{A_1}$. Note that by definition,   $\alpha\in 1+\m_{B_1}$.

Let $K_0/K$ be the maximal unramified subextension of $K_1/K$. Then $L_0=LK_0$ is the maximal unramified subextension of $L$ in $L_1$. Let $A_0$ be the integral closure of $A$ in $K_0$  and let $B_0$ be the integral closure of $A$ in $L_0$. Claim 4 follows from claims 5 -- 7 below.
 \\\\
{\bf Claim 5.} Let $\alpha\in \sT$. Then $A_1[\alpha']$ is stable under the action of $\Gal(L_1/L_0)$. 
\\
{\it Proof of Claim 5:} We take $\alpha'=\gamma (\alpha-1)(\zeta_p-1)^{-1}\in B_1^\times$ such that  $\gamma$ in $A$. For non-trivial $\tau\in \Gal(L_1/L_0)$, we have $\tau(\alpha')= \gamma (\tau(\alpha)-1)(\tau(\zeta_p)-1)^{-1}= \gamma (\alpha^{\tilde\kappa(\tau)}c_{\tau}-1)(\zeta_p^{\kappa(\tau)}-1)^{-1}$. We have that $\zeta_p^{\kappa(\tau)}-1$  belongs to $(\zeta_p-1)\Z_{(p)}[\zeta_p]^\times$, and that $\alpha^{\tilde \kappa(\tau)}c_{\tau}-1\in c_{\tau}-1+(\alpha-1)A_1[\alpha-1]$. We have $c_{\tau}\equiv 1 \bmod (\alpha-1)B_1$ and consequently, $\gamma(c_{\tau}-1)(\zeta_p-1)\in B_1\cap K_1=A_1$. This proves $\tau(\alpha')\in A_1[\alpha']$. 
\\\\
{\bf Claim 6.} Let $\alpha\in \sT$ and take $\alpha'=\gamma (\alpha-1)(\zeta_p-1)^{-1}\in B_1^\times$ such that  $\gamma$ in $A$. Let  $\alpha_0:= N_{L_1/L_0}(\alpha')$. Then
$A_1[\alpha']= A_1[\alpha_0]$. 
\\
In fact, we apply \Cref{lemuni}  by taking $L_1/K_1$ as $L/K$ in the statement and by taking $y=\alpha'$,  $z=\alpha_0$. Then it is sufficient to prove the following claims (6a) and (6b).
\\
{\bf Claim (6a).}  $\sigma(\alpha')(\alpha')^{-1}-1\in \gamma A_1[\alpha']^\times$.
\\
{\bf Claim (6b).} $\sigma(\alpha_0)\alpha_0^{-1}-1\in \gamma A_1[\alpha']^\times$.
\\\\
{\it Proof of (6a):} In fact, $\sigma(\alpha')(\alpha')^{-1}-1 = \sigma(\alpha-1)(\alpha-1)^{-1}-1= (\zeta_p\alpha-1)(\alpha-1)^{-1}-1.$ Since this is equal to $(\zeta_p-1)(\alpha-1)^{-1}\alpha$, we have 6(a).
\\\\
{\it Proof of (6b):} 
Since $A_1[\alpha']$ is finite over  $A_1$ and is an integral domain, it is a local ring. Let $\m_{A_1[\alpha']}$ denote its unique maximal ideal.
Let $n$ be the order of $\Gal(L_1/L_0)$; note that it is invertible in $A$.
\\
Write $(\zeta_p-1)(\alpha-1)^{-1}\alpha=\gamma u$ with $u\in A_1[\alpha']^\times$. We have 
\begin{center}
   $\sigma(\alpha_0)\alpha_0^{-1}-1 = (\prod_{\tau\in \Gal(L_1/L_0)} \tau(\sigma(\alpha')(\alpha')^{-1}))-1 = (\prod_{\tau\in \Gal(L_1/L_0)} (1+ \gamma \tau(u)))-1.$  
\end{center}

\noindent  Since $\tau(u)\in u(1+ \m_{A_1[\alpha']})$, 
we have $\sigma(\alpha_0)\alpha_0^{-1}-1 \in n \gamma u(1+\m_{A_1[\alpha']})$.  This proves (6b). 
\\\\
Note: We recall that $A_1[\alpha']$ is a local ring and $A_1 \subset A_1[\alpha'] \subset B_1$. 
\\
Consequently, $B_0^{\times} \cap A_1[\alpha'] \subset B_1^{\times} \cap A_1[\alpha']=A_1[\alpha']^\times$.
\\
{\bf Claim 7.} Let the notation be as in Claim 6. Then there is $w\in A_0$ such that $\alpha'':=Tr_{L_0/L}(w\alpha_0)$ satisfies $A_0[\alpha_0]= A_0[\alpha'']$.
\\
{\it Proof of Claim 7:}
We apply  \Cref{lemuni}  by taking $L_0/K_0$ as $L/K$ in the statement and by taking $y=\alpha'$, $z=\alpha''$. 
\\
It is then sufficient to prove that there is $w\in A_0$ such that  $(\sigma-1)Tr_{L_0/L}(w \alpha_0)\in \gamma A_0[\alpha']^\times$. 
\\
We have $(\sigma-1)Tr_{L_0/L}(w \alpha_0)= Tr_{L_0/L}(w(\sigma(\alpha_0)-\alpha_0))$, and $\sigma(\alpha_0)-\alpha_0=\gamma w_0$ for some $w_0\in  A[\alpha']^\times\cap L_0=  A_0[\alpha_0]^\times$. There is $w\in A_0$ such that $Tr_{k_0/k}$ ($k_0$ denotes the residue field of $A_0$) sends the residue class of $ww_0$ to a nonzero element of $k$. 
For this $w$, $(\sigma-1)Tr_{L_0/L}(w\alpha_0)$ is an element of $\gamma A_0[\alpha_0]^\times$. 
\\
As in the proof of (6b), we have
 $B^{\times} \cap A_0[\alpha_0] \subset B_0^{\times} \cap A_0[\alpha_0]=A_0[\alpha_0]^\times$, and therefore, we have Claim 7.
\\\\
This concludes the proof of Claim 4.

\subsection{Proof of \Cref{prime} - Part II.}\label{e=l} 
In this subsection, we prove
\subsubsection{}\label{el}  \Cref{prime} in the case $[L:K]=[\Gamma_B:\Gamma_A]$.  
\subsubsection{}
Let $\ell=[L:K]$. Let $\Xi$ be the set of all nonzero elements of  $B$ whose class in $\Gamma_B$ does not belong to $\Gamma_A\subset \Gamma_B$. 
\\\\
 {\bf Claim 8.} If $s\in \Xi$, we have $Bs^{\ell} \subset A[s]$.
\\
{\it Proof of Claim 8:} Take $a\in A$ whose class in $\Gamma_B$ coincides with the class of $s^{\ell}$. That is, $s^{\ell}=ua$ for some $u \in B^{\times}$. 

Let $x\in B$. Then $x=\sum_{i=0}^{\ell-1} x_is^i$ with $x_i\in K$ and $x_is^i\in B$ for  $0\leq i\leq \ell-1$. We have  $x_ia \in B\cap K=A$. Hence $xa = \sum_{i=0}^{\ell-1} (x_ia)s^i \in A[s]$. Applying this to $x=u$, we have $ua=s^{\ell} \in A[s]$. This proves Claim 8.
\subsubsection{}
Next, consider the following condition (C). 
\\
(C) For each $s\in \Xi$, there exists $t\in \Xi$ such that $s\in Bt^{\ell}$.
\\\\
We first assume (C) is satisfied. Take $s\in \Xi$. 
Let $S$ be a finite subset of $B$. Then by the condition (C), there is $t\in \Xi$ such that for all $x=\sum_{i=0}^{\ell-1} x_is^i\in S$ ($x_i\in K$, $x_is^i\in B$, $x_0\in K\cap B=A$) and for $1\leq i \leq \ell-1$, we have 
$Bt^{\ell} \supset Bx_is^i$. By Claim 8, we have $S\subset A[t]$. 
\subsubsection{}
In the rest of this proof of \ref{el}, we assume that (C) is not satisfied. Then there is  $s_1\in \Xi$ such that $Bs_1 \not\subset Bt^{\ell}$ for any $t\in \Xi$. 
Let $\Gamma'_B$ be the subgroup of $\Gamma_B$ consisting of classes of $x\in L^\times$ such that $Bs_1^n \subset Bx \subset Bs_1^{-n}$ for some $n\geq 0$. We have a homomorphism $\lambda: \Gamma_B' \to \R$ to the additive group $\R$ characterized by the following property. Let $x\in L^\times$ and assume that the class $\text{class}(x)$ of $x$ in $\Gamma_B$ belongs to
$\Gamma'_B$. 
\\\\
Let $m,n\in \Z$, $n>0$. Then $Bs_1^m \subset Bx^n$ if $m/n \geq \lambda(x):=\lambda(\text{class}(x))$, and $Bx^n \subset Bs_1^m$ if $m/n \leq \lambda(x)$.  
\\\\
{\bf Claim 9.} Let $\Xi'$ be the subset of $\Xi$ consisting of all elements of whose classes in $\Gamma_B$ belong to $\Gamma_B'$, and let $E$ be the the subgroup of $\R$ generated by $\{\lambda(s')\;|\; s'\in \Xi'\}$. Then $E$ is isomorphic to $\Z$.
\\
{\it Proof of Claim 9:} If $E$ is not isomorphic to $\Z$, $\ell E$ is dense in $\R$. Therefore, there are elements $s_i$ of $\Xi'$ and integers $n_i$ ($2\leq i\leq m$) such that $1 >\sum_{i=2}^m \ell n_i\lambda(s_i) > 1-\ell^{-1}$.  Let $s=s_1\prod_{i=2}^m s_i^{-n_i\ell}$.  Then $0 < \lambda(s) < \ell^{-1}$ and the image of $s$ in $\Gamma_B/\Gamma_A$ is not trivial.  Thus, $s\in \Xi'$ and $Bs^{\ell} \supset Bs_1$, and we reach a contradiction. This proves Claim 9. 
\\\\
Take $a\in L^\times$ whose class in $\Gamma_B$ belongs to $\Gamma_B'\cap \Gamma_A$ and is sent by $\lambda$ to the positive generator of $E$. Let $\bar s \in \Xi'$ and assume that $\lambda (\bar s)$ is $n \lambda(a)$ for some integer $n$. 
Then $(\bar s) - na$ is in the kernel of $\lambda$. Let $s:= a+ (\bar s) - na$.

\noindent Then $s\in \Xi'$ and is sent by $\lambda$ to the positive generator of $E$. 
\\\\
 {\bf Claim 10.} If $s'\in \Xi'$ and $s'|s$ in $B$, then $\lambda(s')=\lambda(s)$. 
\\
{\it Proof of Claim 10:} If $\lambda(s')< \lambda(s)$, we should have $\lambda(s')=0$. Hence, we have $B(s')^{\ell}\supset B s_1$, contradicting our assumption. This proves Claim 10.
\\\\
Let $R=\{a\in A\smallsetminus \{0\}\;|\; a^{-1}s\in B\}$. 
If $a, b\in R$ and if $a|b$, then $A[a^{-1}s]\subset A[b^{-1}s]$. To complete the proof of \ref{el}, it is sufficient to prove $B=\cup_{a\in R} A[a^{-1}s]$. 
\\\\
Let $x\in B$ and write $x=\sum_{i=0}^{\ell-1} x_i s^i$ ($x_i \in K$, $x_is^i\in B$, $x_0\in A$). Let $1\leq i \leq \ell-1$ and assume that $x_i$ does not belong to $A$. Write $x_i= a_i^{-1}$ with $a_i\in A$. 
Note that each $x_is^i$ is an element of $A[a_i^{-1}s^i]$. 
\\
We prove $a_i \in R$. In fact, if $a_i\notin R$, there is  $j$ such that $1\leq j\leq i-1$ and $s^j |a_i|s^{j+1}$ in $B$. Then since $a_is^{-j}, s^{j+1}a_i^{-1}\in \Xi'$ and $a_is^{-j}|s$, $s^{j+1}a_i^{-1}|s$, we have $\lambda(a_is^{-j})=\lambda(s)$ and $\lambda(s^{j+1}a_i^{-1})=\lambda(s)$ by Claim 10. Therefore, $\lambda(a_is^{-j})+\lambda(s^{j+1}a_i^{-1})=\lambda(s)$ should coincide with $2\lambda(s)$, a contradiction. Hence, $a_i\in R$. This proves $x\in A[a^{-1}s]$ for some $a\in R$.

\subsection{Proof of \Cref{prime} - Part III.}\label{elem} We now complete the proof of Theorem \ref{prime}. 
By  \cref{filuni}, \cref{pf*}, and \cref{e=l}, it remains to prove the following case: $[L:K]=[l:k]$. 
\\\\
This is easy to see. In this case, $l$ is generated over $k$ by an element $s$ of $l$. Then $B=A[\tilde s]$ for any lifting $\tilde s$ of $s$ to $B$.

\subsection{Proofs of \Cref{thm2}, \Cref{appl}, and \Cref{dlequiv}.}

\subsubsection{Proof of \Cref{thm2}}\label{pfthm2}  We have $K=L_0 \subset L_1 \subset \dots \subset L_n = L$ 
such that $L_1/K$ is 
unramified and  $L_i/L_{i-1}$ for $2\leq i\leq n$ is a cyclic extension of degree a prime number (see \cite{N} II \S{9}, e.g. 9.12).

For $0\leq i\leq n$, let $B_i$ be the integral closure of $A$ in $L_i$. Then $B_1=A[s_1]$ for some  $s_1\in B_1$ and for $2\leq i\leq n$, $B_i$ is a filtered union of subrings of the form $B_{i-1}[s_i]$ with $s_i\in B_i$ (\Cref{prime}).  From this we see that $B=B_n$ is a filtered union of subrings $B'$ which have the following property:
\\
There is a subring $B'_i$ of $B_i$ for $0\leq i \leq n$ such $A=B_0'\subset B_1'\subset \dots \subset B'_n=B'$ and such that for $1\leq i\leq n$, $B'_i=B'_{i-1}[s_i]$ for some element $s_i$ of $B'_i$ whose monic irreducible polynomial $f_i(T)$ is with coefficients in $B'_{i-1}$. 
\\\\
We have $B_i'\cong B'_{i-1}[T]/(f_i(T))$ and hence $B_i'$ is of complete intersection over $B'_{i-1}$. This shows that $B'$ is of complete intersection over $A$. 
Theorem \ref{thm2} follows from this.

\subsubsection{Proof of \Cref{appl}.}\label{pfappl} Consider $K=L_0\subset L_1 \subset \dots \subset L_n=L$ as in \ref{pfthm2}.  
We may assume that $L/K$ is a power of $p$. 
We proceed by induction on $e(L/K)$. Assume $e(L/K)>1$. 
Then there is $i$ such that $0 \leq  i<n$ and such $e(L_i/K)=1$ 
and $e(L_{i+1}/L_i)>1$. Let $\chi$ be a non-trivial character of 
$\Gal(L_{i+1}/L_i)$ and write $\rsw(\chi)=h^{-1} \otimes w$ with $h$
 a generator of the ideal $\sH$ associated to $L_{i+1}/L_i$ 
 and $w\in k_{B_i} \otimes_{B_i} \Omega^1_{B_i}(\log)$ where $k_{B_i}$ 
 denotes the residue field of $B_i$. Let $\bar w$ be the image of $w$ in 
 $k_{B_i} \otimes_{\Z} \Gamma_{B_i}= k_{B_i} \otimes_{\Z} \Gamma_A$. 
  Write $\bar w=\sum_{j=1}^s b_j \otimes \text{class}(a_j)$ 
 for $b_j\in k_{B_i}$ and $a_j\in K^\times$. Let $A'$ be the Henselization of the integral closure 
 of $A[U_1, \dots, U_s]_{(\m_A)}$ in $K(U_1,\dots, U_s)((a_jU_j)^{1/p} \; (1\leq j\leq s))$ 
 which is a valuation ring. Then (by \ref{types}),  $A'$ is obtained from $A$ by successive log smooth extensions  of type 2. 
 Since $d\log(a_j)=-d\log(U_j)$ in $k_{A'} \otimes_{A'} \Omega^1_{A'}(\log)$,
  the image of $w$ in $k_{A'} \otimes_{\Z} \Gamma_{A'}$ is zero. But the image of $w$ in $k_{A'} \otimes_{A'} \Omega^1_{A'}(\log)$ is not zero. By \ref{dl3}, this shows $e(L_{i+1}K'/L_iK')=1$. Since $e(L_{t+1}K'/L_tK')\leq e(L_{t+1}/L_t)$ for any $0\leq t\leq n-1$  by \ref{124}, we have $e(LK'/K')<e(L/K)$.

\subsubsection{Proof of \Cref{dlequiv}}\label{7.5pf} The equivalence of (ii) and (iii) is in  \Cref{thm2cor}. 
\\\\
We first prove (i) $\Rightarrow$ (ii). By Theorem \ref{appl}, we may assume that $\Gamma_A=\Gamma_B$. Since $L/K$ is defectless, $[L:K]=[k_B:k]$ where $k_B$ denotes the residue field of $B$. Let $(e_i)_i$ be a $k$-basis of $k_B$ and let $(\tilde e_i)_i$ be the lifting of $(e_i)_i$ to $B$. Then $B$ is generated by $\tilde e_i$ as an $A$-module. 
\\\\
We now prove (ii) $\Rightarrow$ (i). Since a log smooth extension is defectless (\ref{lsdl}) in the sense of (4) of \ref{defect}, we may assume that $B$ is a finitely generated $A$-module.
By discussion in  \ref{defect}(1), $L/K$ is defectless.

\section{Properties of upper ramification groups}\label{s:pr}

We prove several important  properties of our upper ramification groups in this section. 
\\
Let $*=\log$ or $\nl$. Let $I$ be a nonzero proper ideal of $\bar A$. 

\begin{thm}\label{LL'2} Let $L/K$ and $L'/K$ be finite Galois extensions such that $L'\subset L$. If the ramification of $L/K$ is $*$-bounded by $I$, then the ramification of $L'/K$ is $*$-bounded by $I$.

\end{thm}

\begin{pf} By the proof of Theorem \ref{thm2}, it is sufficient to prove that if $(S,T)$ and $(S', T')$ are as in the hypothesis of Proposition \ref{Sapro} (5) and if $(S,T, I)$ separates $L/K$, 
then $(S', T', I)$ separates $L'/K$. This follows from Proposition \ref{Sapro} (5). 
 \end{pf}

\begin{pr}\label{M1}  Let $M\subset \bar K$ be the union of all finite Galois extensions $M'/K$ in $\bar K$ such that the ramification of each  $M'/K$ is $*$-bounded by $I$. Then $\Gal(\bar K/K)_*^I= \Gal(\bar K/M)$. 

\end{pr}

\begin{pf} Use \ref{LL'} and Theorem \ref{LL'2}. \end{pf}

\begin{pr}\label{M2} Let $L$ be a finite Galois extension of $K$. Then 
$\Gal(L/K)_*^I=\{1\}$  if and only if the ramification of $L/K$ is $*$-bounded by $I$.

\end{pr}

\begin{pf} This follows from Theorem \ref{M1}.
\end{pf}

\begin{pr}\label{st} Let $s$ be an element of $B$ such that $L=K(s)$, let $f$ be the monic polynomial of $s$ over $K$, and let $S=\{s\}$, $T=\{t\}$ where $t=f(y_s)$. Let $J$ be the ideal of $\bar A$ consisting of all $a\in \bar A$ such that $|\phi(s)-\phi'(s)|_{\bar A}>|a|_{\bar A}$ if $\phi, \phi'\in \Phi(L/K)$ and $\phi\neq \phi'$. Let $I$ be the ideal of $\bar A$ generated by 
$a^{[L:K]}$ for all $a\in J$. Then $(S,T, I)$ separates $L/K$. 
\end{pr}

\begin{pf} For $\phi\in \Phi(L/K)$, let 
 $$U_{\phi}:= \{z\in X(S, T, I)\;|\; |y_s- \phi(s)|(z) \leq |a|_{\bar A}\;\text{for some}\;a\in J\}.$$ 
 Then $U_{\phi}$ is an open subset of $X(S, T,I)$ and $\phi\in U_{\phi}$. We will now show that $X(S, T, I)$ is the disjoint union of $U_{\phi}$ for $\phi\in \Phi(L/K)$.\\ 
 Note that
$t= \prod_{\phi\in \Phi(L/K)} (y_s-\phi(s))$. If $z\in X(S, T, I)$ and if $z\notin U_{\phi}$ for any $\phi\in \Phi(L/K)$, then since $|y_s-\phi(s)|(z) > |a|_{\bar A}$ for all $\phi\in \Phi(L/K)$ and all $a\in J$, we have $|t|(z)> |a|_{\bar A}^{[L:K]}$ for any $a\in J$ and this contradicts $z\in X(S, T, I)$. This shows $X(S, T, I)=\cup_{\phi}\;U_{\phi}$.
\\  
If $\phi'\in \Phi(L/K)$ and  $z\in U_{\phi}\cap U_{\phi'}$, then $|y_s-\phi(s)|(z)\leq |a|_{\bar A}$ and $|y_s-\phi'(s)|(z)\leq |a|_{\bar A}$ for some $a\in J$, and we have  $|\phi(s)-\phi'(s)|(z) \leq |a|_{\bar A}$, hence $\phi =\phi'$.
\end{pf}

\begin{pr} Let $L/K$ be a finite Galois extension. For $*=\log, \nl$, 
$\Gal(L/K)^I_*=\{1\}$ if $I$ is sufficiently small.

\end{pr}

\begin{pf} For $*=\nl$, this follows from the propositions \ref{M2} and \ref{st}. The case $*=\log$ follows from the case $*=\nl$ by Proposition \ref{easy}. 
\end{pf}

\begin{pr}\label{lnl} Let $L/K$ be a finite Galois extension and let $A'/A$ be a log smooth extension of Henselian valuation rings with field of fractions $K'$ such that  $e(LK'/K')=1$ (Theorem \ref{appl}). Then the restriction $$\Gal(LK'/K')\to \Gal(L/K)$$ induces an isomorphism 
$$\Gal(LK'/K')_{\nl}^I\overset{\cong}\to \Gal(L/K)_{\log}^I$$
 for any nonzero proper ideal $I$ of $\bar A$.
\end{pr}
\begin{pf}
This follows from \ref{124}. 
\end{pf}

\begin{thm}\label{cap} Let $\Lambda$ be a non-empty set of nonzero proper ideals of $\bar A$ and let $J=\cap_{I\in \Lambda} \;I$. Then $\Gal(\bar K/K)_*^J=  \cap_{I\in \Lambda} \Gal(\bar K/K)_*^I$. 
\end{thm} 

\begin{pf} Assume  $*=\nl$. Let $L/K$ be a finite Galois extension and let $(S,T)$ be a pre-presentation for $L/K$ which is of complete intersection. 
 It is sufficient to prove that if $(S, T, J)$ separates $L/K$, then for some $I\in \La$, $(S,T, I)$ separates $L/K$. By Proposition \ref{Sapro} (4), there is an ideal $I_1$ of $\bar A$ such that for each nonzero proper ideal $I$ of $\bar A$, $(S, T, I)$ does not separate $L/K$ if and only if  $I\supset I_1$. If $(S, T, I)$ does not separate $L/K$ for every $I\in \La$, then $I\supset I_1$ for every $I\in \La$. Consequently,  $J=\cap_{I\in \La}\; I \supset  I_1$, and this implies that $(S, T, J)$ does not separate $L/K$.
 
 The log version follows from the non-log version. 
\end{pf}

\begin{pr}\label{Ii} Let $L$ be a finite Galois extension of $K$. Then there is a finite sequence of ideals 
$0 \subsetneq  I_1\subsetneq \dots \subsetneq I_n\subsetneq \bar A$ of nonzero proper ideals of $\bar A$ such that for any nonzero proper ideal $I$ of $\bar A$, we have:

$$\Gal(L/K)_*^I=\{1\}\te{if} I\subsetneq I_1,$$

$$\Gal(L/K)_*^I=\Gal(L/K)_*^{I_i} \te{if} 1\leq i <n \te{and} I_i\subset I\subsetneq I_{i+1},$$ 

$$\Gal(L/K)_*^I=\Gal(L/K)_*^{I_n} \te{if} I_n\subset I.$$

\end{pr}

\begin{pf} Let $\{1\} =H_0\subsetneq H_1\subsetneq \dots \subsetneq H_n=\Gal(\bar K/K)$ be the set of all subgroups of $\Gal(L/K)$ which are equal to  $\Gal(L/K)_*^I$ for some nonzero proper deal $I$ of $\bar A$. 
\\
For each $1\leq i\leq n$, let $I_i$ be the intersection of all nonzero proper ideals $I$ of $\bar A$ such that $\Gal(L/K)_*^I=H_i$. By Theorem \ref{cap}, these $I_i$'s have the required properties. 
\end{pf}

\noindent We will prove the following proposition in section \ref{s:th1}.
\begin{pr}\label{mA}
For $I=\m_{\bar A}$, $\Gal(\bar K/K)_{\log}^I$ is the wild inertia group, and $\Gal(\bar K/K)^I_{\nl}$ is the inertia group.

\end{pr}

\section{Theory of Abbes-Saito}\label{s:AS}

In this section,  we will prove the relation of our upper ramification groups and those   of Abbes-Saito stated in \ref{relAS}.
\\\\
For this entire section,  we assume that $A$ is a complete discrete valuation ring. 
\\
In \ref{XtS} and \ref{XtS2}, let $r\in \Q_{>0}$ and let  $I= I(r)=\{x\in \bar A\;|\;\text{ord}_{\bar A}(x) \geq r\}$.
\\In \ref{relAS8} we expand our discussion to include non-principal ideals of $\bar A$.
\subsection{Ramification bounded by $I=I(r)$.}\label{XtS} 

Let $L/K$ be a finite Galois extension and let $B$ be the integral closure of $A$ in $L$.

Let $(S, T)$ be a presentation of $B/A$ (\ref{pres}). Then as a topological space, 
$X(S, T, I)$ coincides with the Zariski-Riemann space associated with $X^r_S$ in \cite{AS1} used by Abbes-Saito in their definition of the non-log upper ramification groups (\ref{adican}). 
 In particular, the ramification of $L/K$ is non-logarithmically bounded by $I$ in our sense if and only if it is bounded by $r$ in the  sense of Abbes-Saito \cite{AS1}.

By Theorem \ref{appl} and by \ref{124}, the ramification of $L/K$ is logarithmically bounded by $I$ if and only if 
 there is a sequence of valuation rings $A=A_0\subset A_1\subset \dots \subset A_n=A'$ 
 satisfying the following:

 \begin{itemize}
     \item Each $A_i/A_{i-1}$ for $i=1$ (resp. $2\leq i\leq n$) is a log smooth extension of Henselian valuation rings of type 1 (resp. type 2).
     \item The ramification index $e(LK'/K')=1$,  where $K'$ is the field of fractions of $A'$ and such that $LK'/K'$ has ramification non-logarithmically bounded by $I$. 
     
 \end{itemize}
 
 \noindent In this case $A'$ is a discrete valuation ring. 
  By \cite[1.2.6]{Sa0}, this shows that the ramification of $L/K$ is logarithmically bounded by $I$ in our sense if and only if it is logarithmically bounded by $r$ in the sense of Abbes-Saito \cite{AS1}.

\subsection{Comparison of the filtrations.}\label{XtS2} Let $G=\Gal(\bar K/K)$.
\\
Recall that 
 the definitions of $G^r_{\log}$ and $G^r$ by Abbes-Saito in \cite{AS1} are as follows.\\
 The group $G^r_{\log}$ (resp. $G^r$) is the intersection of kernels of $G\to \Gal(L/K)$ where $L$ ranges over all finite Galois extensions $L$ of $K$ in $\bar K$ such that the ramification of $L/K$ is non-logarithmically (resp. logarithmically) bounded by $r$  in the sense of Abbes-Saito. 
 \\\\
 Hence, by \ref{M2} and \ref{XtS}, our $G^I_{\log}$ (resp. $G^I_{\nl}$) for $I=I(r)$ coincides with their $G^r_{\log}$ (resp. $G^r$).

\subsection{$G^I_*$ for general $I$.}\label{relAS8} Let $I$ be a nonzero proper  ideal of $\bar A$ which need not be principal.  
\\
By Proposition \ref{Sapro} (2) applied to a presentation $(S,T)$ for $L/K$, $G^I_{\nl}$ 
 coincides with the closure of $\cup_J\; G^J_{\nl}$ in $G$ where $J$ ranges over all principal nonzero subideals of $I$. 
 From this, we also have  an analogous  coincidence for the logarithmic filtration. 
That is, for $*= \nl \te{or} \log$, 
\\\\$G^I_{*}$ 
 coincides with the closure of $\cup_J\; G^J_{*}$ in $G$ where $J$ ranges over all principal nonzero subideals of $I$.

\section{Proof of \Cref{thm1}}\label{s:th1}

 In this section, we recall and prove  \Cref{thm1}, and we also give the proof of \Cref{mA} in \ref{pfma}, as an application of \Cref{thm1}. 
 \\ Except in \ref{pfma}, let the assumptions be as in Theorem \ref{thm1}; we recall its statement below for convenience.

\begin{Thm}[\ref{thm1}] Assume that the residue field of $A$ is of characteristic $p>0$. Let $L$ be a cyclic extension of $K$ of degree $p$ and let $\sH\subset A$ be the associated ideal. Then for a nonzero proper ideal $I$ of $\bar A$, the image of $G_{\log}^I$ in $\Gal(L/K)$ is $\Gal(L/K)$ if and only if $I\cap A\supset \sH$ and is $\{1\}$ if and only if $I\cap A \subsetneq \sH$.
\end{Thm}

\subsection{Preparation.}

\begin{sbprop}\label{notsep1} Let $s\in B\smallsetminus A$. Let $S=\{s\}$, $T=\{t\}$ with $t= \prod_{\phi\in \Phi(L/K)} (y_s-\phi(s))$ in $A[y_s]$. Taking $\phi_0,\phi_1\in \Phi(L/K)$, $\phi_0\neq \phi_1$, let $b=\phi_1(s)-\phi_0(s)$ (note that $|b|_{\bar A}$ is independent of the choices of $\phi_0, \phi_1$), and let $I=\bar A b^p$. Then
 the space $X(S, T, I)$ is connected.  In particular, $(S, T, I)$ does not separate $L/K$. 

\end{sbprop}

\begin{pf} 

The space $X(S, T, I)$ is the subspace of $D^S$ consisting of all $z\in D^S$ such that
\begin{center}
   $\prod_{\phi\in \Phi(L/K)} |y_s- \phi(s)|(z)\leq |b|^p_{\bar A}$ (see \cref{XSTI}). 
\end{center}

 Let $y=y_s-\phi_0(s)$. Then $X(S, T, I)$ is the set of all $z\in D^S$ such that \begin{center}$|y \prod_{\phi\in \Phi(L/K)\smallsetminus \{\phi_0\}} (y+\phi_0(s)-\phi(s))|(z)\leq |b|^p_{\bar A}$. \end{center} 

\noindent {\bf Claim 11.} $X(S, T, I)$ is the set of all $z\in D^S$ such that  $|y|(z) \leq |b|_{\bar A}$.\\
{\it Proof of Claim 11:} If $|y|(z)\leq |b|_{\bar A}$, then for any $\phi\in \Phi(L/K)\smallsetminus \{\phi_0\}$, $|\phi_0(s)-\phi(s)|(z)=|b|_{\bar A}$ and hence, $|y+\phi_0(s)-\phi(s)|(z)\leq |b|_{\bar A}$. As a consequence, $|y\prod_{\phi\in \Phi(L/K)\smallsetminus \{\phi_0\}} (y+\phi_0(s)-\phi(s))|(z)\leq |b|^p_{\bar A}$. Assume $|y|(z)>|b|_{\bar A}$. Then $|y+\phi_0(s)-\phi(s)|(z)>|b|_{\bar A}$ and hence, $|y \prod_{\phi\in \Phi(L/K)\smallsetminus \{\phi_0\}} (y+\phi_0(s)-\phi(s))|(z)>|b|^p_{\bar A}$.  
\\\\
By Claim 11, $X(S, T, I)$ is connected.
\end{pf}

\begin{sbprop}\label{sep1}  In the situation of Proposition  \ref{notsep1},  $(S, T, \m_{\bar A}b^p)$ separates $L/K$.

\end{sbprop}

\begin{pf} This follows from Proposition \ref{st}. 
\end{pf}

\subsection{The proof.}\label{pfth1} 

We now prove \Cref{thm1}.

\begin{proof}
Let the notation be as in \ref{notsep1} and we will denote (a choice of) $b$ corresponding to any fixed $s \in B \backslash A$ by $b_s$. (As the following argument only depends on $|b_s|_{\bar{A}}$, this choice does not matter.)
First we will compare $\bar{A}\sH$ and $\bar{A}b_s^p$ when $e(L/K)=1$. 
\\

The ideal $\sH$ of $A$ is generated by  $\{N_{L/K}(x)\;|\;x\in \js\}$. As in the proof of \ref{HH'}, we recall that when $e(L/K)=1$, the ideals $\is$ and $\js$ of $L/K$ coincide. (Such basic properties of the ideals $\js, \is$ of $B$ are explored in \cite{V1, V2}.) \\ 

Consequently,  $\sH$  is the ideal of $A$ generated by  $\{N_{L/K}(x)\;|\;x\in \is\}$. The ideal $\is$ is generated by $\{b_s \}_s$, by definition. Therefore, in this case, $\bar{A} \sH = \bigcup_{s \in B \backslash A} \bar{A}b_s^p$.  
\\

In the defectless case, by Theorem \ref{appl} and Proposition \ref{HH'}, we may assume that $B=A[s]$ for some $s\in B$ and $e(L/K)=1$. 
Then Theorem \ref{thm1} in this case follows from the propositions \ref{lnl}, \ref{notsep1} and \ref{sep1}. We note that in this case, $\bar{A} \sH =  \bar{A}b_s^p$ for this particular $s$.

In the defect case, $e(L/K)=1=f(L/K)$ by definition. Then Theorem \ref{thm1} follows from Theorem \ref{prime} and the propositions \ref{HH'}, \ref{notsep1}, and \ref{sep1}. \end{proof}

\subsection{\Cref{mA} as an application of \Cref{thm1}.}\label{pfma}

We want to prove that $\Gal(\bar K/K)_{\log}^{\m_{\bar A}}$ is the wild inertia group, and $\Gal(\bar K/K)^{\m_{\bar A}}_{\nl}$ is the inertia group.

\begin{pf}
The non-log case is easy. 
We prove the log case.

 It is clear that the ramification of a finite tame Galois extension of $K$  is logarithmically bounded by $\m_{\bar A}$.

We prove that the ramification of a non-tame finite Galois extension $L/K$  is not logarithmically bounded by $\m_{\bar A}$. 

Let $p$ be the characteristic of the residue field $k$. Replacing $K$ by its finite tame extension, we may assume that we  have a surjection $\Gal(L/K)\to \Gal(L'/K)=\Z/p\Z$; $L'\subset L$ and $L'/K$ is not unramified.  Then the ideal $\sH$ of $A$ associated to $L'/K$ is contained in  $\m_A$. By Theorem \ref{thm1}, the ramification of $L'/K$ is not logarithmically bounded by  $\bar A \sH$. Hence, the ramification of $L'/K$ is not logarithmically bounded by $\m_{\bar A}$. 
\end{pf}

\section{Breaks of the (log) filtration}\label{s:br}
In this section, we consider logarithmic upper ramification groups. 
\subsection{Definitions.}
\begin{sbdefn}
  For a nonzero proper ideal $I$ of $\bar A$, we say {\it $I$ is a break (of the logarithmic upper ramification filtration) of $A$}  if  $\Gal(\bar K/K)^I_{\log}$ does not coincide with the closure of the union of  $\Gal(\bar K/K)^J_{\log}$ for all nonzero ideal $J\subsetneq I$.  
\end{sbdefn} 
\begin{sbdefn}
  For a finite Galois extension $L/K$, we say {\it $I$ is a  break for $L/K$} if for every nonzero ideal $J\subsetneq I$ of $\bar A$, $\Gal(L/K)^I_{\log}\neq \Gal(L/K)^J_{\log}$. 
 
\end{sbdefn}
\noindent Thus $I$ is a break of $A$ if and only if $I$ is a break of some finite Galois extension $L/K$.

\subsection{Breaks for defectless extensions.}
Abbes-Saito (\cite{AS1}) proved that in the case $A$ is a discrete valuation ring, if $I$ is a break of the upper ramification filtration,  then $I$ is a principal ideal. 
We generalize this to the following theorem.

\begin{thm}\label{Idl2} Let $L/K$ be a defectless finite Galois extension. Then any  break of  $L/K$ is a principal ideal. 

\end{thm}

\begin{pf} By  Theorem \ref{dlequiv}, this follows from 
   Proposition \ref{Sapro} (4).
  \end{pf}

\subsection{Breaks for defect extensions.}
 In this subsection, we show that various types of breaks appear when we allow defect. 
 \\
{\bf Note:} Only in this subsection, the valuation of a valuation ring is treated additively, not multiplicatively. This is 
because an important case of this section is that the value group is a subgroup of the additive group $\R$. That is, if $v_A(a)$ denotes the valuation of $a \in K$, then   $v_A(ab)=v_A(a)+v_A(b)$ for all $a, b \in K$ and for $a\in K$,  $b\in K^\times$, $v_A(a)\geq v_A(b)$  means $b^{-1}a\in A$.
\subsubsection{Set-up.}\label{se1}
For a valuation ring $A$ with value group $\Gamma$ (written additively), there is a bijection 
$$I\mapsto C:= \{v_A(x)\;|\; x\in I\smallsetminus \{0\}\}$$
from the set of  ideals of $A$ to the set of all subsets $C$ of $\Gamma_{\geq 0}:= \{\gamma \in \Gamma\;|\;\gamma \geq 0\}$ satisfying the following condition (i). 

\medskip

(i) If $\gamma\in  C$, any element  $\gamma'$ of $\Gamma$ such that  $\gamma'\geq \gamma$ belongs to $C$.

\begin{sbprop}\label{se2} Let $\Gamma$ be a totally ordered abelian group (written additively). Let $C$ be a subset of $\Gamma_{\geq 0}$ satisfying the condition (i) in \ref{se1} such that $C\neq \emptyset, \Gamma_{\geq 0}$. Let $p$ be a prime number. Then the following (a) and (b) are equivalent.

\medskip

\noindent (a)  Either $C$ has a minimal element or for each $c\in C$, there is $d \in \Gamma$ such that $pd\in C$ and $c\geq pd$. 

\medskip

\noindent (b) There are a Henselian valuation ring $A$ of characteristic $p$ whose value group is $\Gamma$ and a cyclic  extension $L$ of the field of fractions $K$ of $A$ of degree $p$ satisfying the following condition.

The set $C$ corresponds to the ideal $\sH$ of $A$ associated to $L/K$ (\ref{rev1}). That is (Theorem \ref{thm1}), if $I$ is a nonzero proper ideal of $\bar A$, $\Gal(L/K)_{\log}^I=\Gal(L/K)$ if and only if the subset $C'$ of $\Gamma_{\geq 0}$ corresponding to the ideal $I\cap A$ of $A$ satisfies $C'\supset C$.

\end{sbprop}

\noindent This is proved later in \ref{se8}, \ref{pf31}, \ref{pf32}, and \ref{pf33}.

\begin{sbrem}
(1) Let $\Gamma= \Z^2$ with the lexicographic order. Then for a prime number $p$, the set $C=\{(m,n)\in \Z^2\; |\; m>0\}$ satisfies the condition (i) in \ref{se1} but does not satisfy the condition $(a)$ in \ref{se2}. 
\\
(2) If $\Gamma$ is a nonzero subgroup of $\R$, any subset $C$ of $\Gamma_{\geq 0}$ satisfying (i) in \ref{se1} satisfies the condition $(a)$. This is because if $\Gamma$ is not isomorphic to $\Z$, $p\Gamma$ is dense in $\R$.\\\\ 
Thus, Proposition  \ref{se2} shows the following propositions \ref{se4} and \ref{se7}. 

\end{sbrem}

\begin{sbprop}\label{se4}  Let $\Gamma$ be a  nonzero subgroup of $\R$ which is not isomorphic to $\Z$. Let $p$ be a prime number.  Let $a\in \R_{> 0}$ and assume $a\in \Gamma$ (resp. $a\in \Gamma$, resp. $a\notin \Gamma$). Then  there are a Henselian valuation ring $A$ of characteristic $p$ whose value group is $\Gamma$ and a cyclic  extension $L$ of the  field of fractions $K$ of $A$ of degree $p$ such that for each nonzero proper ideal $I$ of $\bar A$, $\Gal(L/K)_{\log}^I=\Gal(L/K)$ if and only if the subset $C'$ of $\Gamma$ corresponding to $I\cap A$ satisfies $C'\supset C$ where $C=\{x\in \Gamma\; |\; x\geq a\}$ (resp. $C=\{x\in \Gamma\;|\; x>a\}$, resp. $C=\{x\in \Gamma\;|\; x>a\}$). 

\end{sbprop}

\begin{sbprop}\label{se7}  Let $\Gamma$ be a  nonzero subgroup of $\R$ which is not isomorphic to $\Z$.  Let $a\in \R_{> 0}$. Then  there is a Henselian valuation ring $A$ whose value group is $\Gamma$ such that 
the ideal $\{x\in \bar A\;|\; v_{\bar A}(x)>a\}$ of $\bar A$ is a break of $A$. 

\end{sbprop}

\subsubsection{Preparation for the proof of \Cref{se2}.}\label{se8} 
In general, let $\Gamma$ be a totally ordered abelian group whose group law is written additively, let $R_0$ be an integral domain, and let  $R$ be the group ring of $\Gamma$ over $R_0$. We will denote the group element of $R$ corresponding to $\gamma\in \Gamma$ by $t^{\gamma}$. It follows that
$t^{\gamma}t^{\gamma'}= t^{\gamma+\gamma'}$ for $\gamma, \gamma'\in \Gamma$. For a nonzero finite sum $b=\sum_{\gamma \in \Gamma} a_{\gamma}t^{\gamma}\in R$  ($a_\gamma\in R_0$), let $v(b):= \inf\{\gamma\;|\;a_{\gamma}\neq 0\}\in \Gamma$.  Then $v(bb')=v(b)+v(b')$ for nonzero elements $b, b'\in 
R$. Let $Q(R_0)$ and $Q(R)$ be the  fields of fractions of the integral domains $R_0$ and $R$, respectively. Then $v$ extends to a valuation of $Q(R)$. Let $V\subset Q(R)$ be the valuation ring of $v$. Then the value group of $V$ is identified with $\Gamma$ and the residue field of $V$ is identified with $Q(R_0)$.

\subsubsection{Proof of \ref{se2} $(a)\Rightarrow (b)$ assuming that $C$ has a minimal element $c$.}\label{pf31}
Assume first that $c$ does not belong to $p\Gamma$. Take a Henselian valuation ring $A$ of characteristic $p$ with value group $\Gamma$. Let $h$ be an element of $A$ such that $v_A(h)=c$ and let $L=K(\alpha)$,  $\alpha^p-\alpha=1/h$. Then the ideal $\sH$ associated to $L/K$ corresponds to $C$. 
\\\\
Assume next $c=pd$ for some $d\in \Gamma$. Take a Henselian valuation ring $A$ of characteristic $p$ with value group $\Gamma$ and with imperfect residue field. Let $h$ be an element of $A$ such that $v_A(h)=d$, let $u$ be an element of $A$ whose residue class is not a $p$-th power, and let $L=K(\alpha)$, $\alpha^p-\alpha=u/h^p$. Then the ideal $\sH$ associated to  $L/K$ corresponds to $C$.

\subsubsection{Proof of \ref{se2} $(a)\Rightarrow (b)$ assuming that $C$ has no minimal element.}\label{pf32} 
Let $D=\{\gamma \in \Gamma\;|\; p\gamma \in C\}$ and  let $R_0$ be the polynomial ring over $\F_p$ in variables $x_d$ ($d\in D$). Consider the rings $R$ in \ref{se8} and the valuation $v$ on $R$, and let $L':=Q(R)$ and let $B':=V$ there. Then the value group of $B'$ is $\Gamma$ and the residue field of $B'$ is the pure transcendental extension of $\F_p$ with transcendence basis $x_d$ ($d\in D$). The following statement (Claim 12) is proved easily.
\\{\bf Claim 12.} $D$ satisfies the condition (i) in \ref{se1} (when we replace $C$ by $D$ in this condition) and has no minimal element.

\medskip
Let $\sigma: R \to R$ be the ring homomorphism defined by $\sigma(x_d)= x_d+t^d$ for $d\in D$ and $\sigma(t^{\gamma})= t^{\gamma}$ for all $\gamma\in \Gamma$. Then this is an automorphism of $R$ which does not change $v$. Hence, it  induces automorphisms of $L'$ and $B'$. We have $\sigma^p=1$. Let $K'$ be the $\sigma$-fixed part of $L'$ and let $A'$ be the $\sigma$-fixed part of $B'$.  Then $A'$ is a valuation ring and $K'$ is the field of fractions of $A'$. Since $t^{\gamma}\in K'$ for all $\gamma\in \Gamma$, the value group of $A'$ is $\Gamma$. We show that the residue field of $A'$ is the same as that of $B'$. For $d\in D$, take $d'\in D$ such that $d>d'$ (Claim 12).  Then since $\sigma(x_dt^{-d})=x_dt^{-d}+1$ and $\sigma(x_{d'}t^{-d'})= x_{d'}t^{-d'}+1$, we have $x_dt^{-d}-x_{d'}t^{-d'}\in K'$ and therefore, $y_d:=x_d- x_{d'}t^{d-d'}\in K'$. This $y_d$ belongs to $A'$ and the residue class of $y_d$ coincides with that of $x_d$. Hence, the residue field of $B'$ coincides with that of $A'$.

Let $A$ be the Henselization of $A'$, let $K$ be the field of fractions of $A$, and let $B=B'\otimes_{A'} A$. Then $B$ is the integral closure of $A$ in the field $L=L'\otimes_{K'} K$ and is a valuation ring. 
\\
Now we consider the ideal $\js$ of $B$ (\ref{rev1}). 
\\\\
{\bf Claim 13.} $\js$ coincides with the ideal of $B$ generated by $t^d$ for $d\in D$.\\
{\it Proof of Claim 13:}  Since $A$ and $B$ have the same  value group,  $\js$ coincides with ideal generated by $\sigma(f)-f$ $(f\in B)$. 
\\
Clearly, we have:
\\
(1) The ideal $(\sigma(f)-f\;|\; f\in R)$ of $R$ is generated by $t^d$ ($d\in D$). 
\\
Furthermore, as is easily seen, we have:
\\
(2) For a nonzero element $f$ of $R$, $\sigma(f)-f$ belongs to the ideal of $B'$ generated by $t^{v(f)+d}$ ($d\in D$). 
\\
For $f, g\in R$ such that $g\neq 0$ and $v(f)\geq v(g)$, since $N(g)g^{-1}\in R$ ($N$ denotes the norm map of $L'/K'$) and $fg^{-1}= (f\cdot N(g)g^{-1})N(g)^{-1}$, we have:
\\
(3) $B'$ coincides the set of elements of the form $fg^{-1}$ such that $f, g\in R$, $g\neq 0$, $\sigma(g)=g$,  $v(f)\geq v(g)$.
\\
\\
Claim 13 follows from (1) -- (3). 
\\\\
By \ref{rev1}, the ideal $\sH$ of $A$ is generated by $t^{pd}$ ($d\in D$). This completes the proof of $(a) \Rightarrow  (b)$.

\subsubsection{Proof of \ref{se2} $(b)\Rightarrow (a)$}\label{pf33} 

We may assume that $C$ has no minimal element. 
\\
As in \ref{rev1}, $\sH$ is generated by norms of elements of $\js$. So if $D$ denotes the subset of $\Gamma_{\geq 0}$ corresponding to $\js$, we have $C= \{\gamma \in \Gamma\; |\; \gamma \geq pd \; \text{for some}\; d\in D\}$. Hence $C$ satisfies $(a)$. 
\\\\
This completes the proof of \Cref{se2}.\\\\ Now we will present the following improvements of some of the above results under certain assumptions.
\begin{sbprop}\label{se5} In the hypothesis of Proposition \ref{se2} (resp. Proposition \ref{se4}), assume that the dimension $n$ of $\Gamma \otimes \Q$ over $\Q$ is finite. Then we can improve \ref{se2} (resp. Proposition \ref{se4}) by adding  the following  condition (**) to the condition $(b)$ in \ref{se2} (resp. to the conditions on $A$). 

\medskip

(**) There is an algebraically closed subfield $k$ of $A$ such that the field of fractions of $A$ is of transcendence degree $\leq n+1$ over $k$.

\end{sbprop}

\begin{pf} Assume first that $C$ has a minimal element $c$.
\\
If $c$ does not belong to $p\Gamma$, then in \ref{se8}, let $R_0$ be an algebraically closed field $k$ of characteristic $p$, let $A$ in \ref{pf31} be the Henselization of $V$, and let $h$ in \ref{pf31} be $t^c$. 
\\
If $c=pd$ with $d\in \Gamma$, then in \ref{se8}, let $R_0=k(U)$ with $k$ an algebraically closed field of characteristic $p$ and with $U$ an indeterminate, let $A$ in \ref{pf31} be the Henselization of $V$,  let $h$ in \ref{pf31} be $t^d$, and let $u$ in \ref{pf31} be $U$. 
Then the conditions $(b)$ and  (**) are satisfied.

Assume next that $C$ has no minimal element.
\\
In \ref{pf32}, take $y_d$ ($d\in D$) more carefully as follows. Since $\dim_{\Q}(\Gamma\otimes \Q)$ is finite, there are elements  $e(i)$ ($i=0,1,2,\dots$) of $D$ such that $e(0)>e(1)>e(2)>\dots$ and such that for each $d\in D$, there is $i$ such that $d>e(i)$. For $i\geq 0$, define $y_{e(i)}= x_{e(i)}- x_{e(i+1)}t^{e(i)-e(i+1)}$. For $d\in D$ which does not belong to $\{e(i)\;|\;i\geq 0\}$, take $i\geq 0$ such that $d>e(i)$ and let $y_d=x_d-x_{e(i)}t^{d-e(i)}$. 
Let $k'= \F_p(y_d\;;\;d\in D)$ and let $k$ be an algebraic closure of $k'$. Let $A_1$,  $B_1$, $K_1$, $L_1$  be the $A$, $B$, $K$, $L$ in \ref{pf32}, respectively,  and let $A_2= A_1\otimes_{k'} k$, $B_2=B_1\otimes_{k'} k$, $K_2=K_1\otimes_{k'} k$, $L_2=L_1\otimes_{k'} k$. Then $A_2$ and $B_2$ are Henselian valuation rings. If we use $A_2$ and the extension $L_2/K_2$ as $A$ and $L/K$, $(b)$ and (**) are satisfied. In fact, 
 let $z= (x_{e(0)}t^{-e(0)})^p- x_{e(0)}t^{-e(0)}\in K$ and let $\gamma_i$ ($1\leq i\leq n$) be elements of $\Gamma$ which form a $\Q$-basis of $\Gamma \otimes \Q$. Then $L_2$ is algebraic over $k(z, t^{\gamma_1}, \dots, t^{\gamma_n})$ and hence $K_2$ is also algebraic over $k(z, t^{\gamma_1}, \dots, t^{\gamma_n})$. 
 \\\\
 This proves the part of 
\ref{se5} concerning \ref{se2}. 
\\\\By our proof of Proposition \ref{se4} using \ref{se2}, this improvement of \ref{se2} gives the desired improvement of Proposition \ref{se4}. 
\end{pf}

\subsection{Concerning \ref{why}.}\label{br10} Here we explain that $A$ with $L_1/K$ and $L_2/K$ as in \ref{why} exists. And hence, we indeed need to consider all nonzero proper ideals of $\bar A$ in our indexing set for the filtration.
\\\\
Take a Henselian valuation ring $A$ of  characteristic $p$ whose value group is a nonzero subgroup of $\R$ which is not isomorphic to $\Z$ and whose residue field is not perfect,  and an Artin-Schreier extension $L_2$ of $K$ of degree $p$ such that the ideal $\sH$ of $A$ associated to $L_2/K$ is $b^p{\mathfrak m}_A$ for some nonzero element $b$ of  ${\mathfrak m}_A$. Such $A$ and $L_2/K$ exists by Proposition  \ref{se4}. Take a unit $u$ of $A$ whose residue class  is not a $p$-th power and let $L_1= K(\beta)$ where $\beta$ is a solution of $\beta^p-\beta= ub^{-p}$. 
Then the ideal $\sH$ of $A$ associated to $L_1/K$ is $b^pA$. This $(A, L_1/K, L_2/K)$ satisfies the condition in \ref{why}.

\subsection*{Acknowledgments.} 
We are grateful to Takeshi Saito  for his comments on an earlier draft of the manuscript.
\\
K. Kato is partially supported by NSF Award 1601861.

\medskip
\medskip
\medskip \noindent Kazuya Kato\\Department of Mathematics\\University of Chicago\\Chicago, Illinois 60637\\Email: \texttt{kkato@math.uchicago.edu}

\medskip
\medskip
\medskip \noindent Vaidehee Thatte\\Department of Mathematics\\King's College London\\Strand, London WC2R 2LS, United Kingdom\\
Email: \texttt{vaidehee.thatte@kcl.ac.uk}

\end{document}